\documentclass[a4paper]{article}
\usepackage{amsmath,amsthm,amssymb,amscd}

\numberwithin{equation}{section}

{\theoremstyle{definition}\newtheorem{definition}{Definition}[section]
\newtheorem{notation}[definition]{Notation}
\newtheorem{terminology}[definition]{Terminology}
\newtheorem{remark}[definition]{Remark}
}
\newtheorem{proposition}[definition]{Proposition}
\newtheorem{lemma}[definition]{Lemma}
\newtheorem{theorem}[definition]{Theorem}
\newtheorem{corollary}[definition]{Corollary}

\newcommand{\cG}{\mathbb{G}}
\newcommand{\de}{\Delta}
\newcommand{\ot}{\otimes}
\newcommand{\cGh}{\widehat{\mathbb{G}}}
\newcommand{\recht}{\rightarrow}
\newcommand{\B}{\operatorname{B}}
\newcommand{\oH}{\overline{H}}
\newcommand{\oK}{\overline{K}}
\newcommand{\sde}{\delta}
\newcommand{\vphi}{\varphi}
\newcommand{\Mor}{\operatorname{Mor}}
\newcommand{\Hom}{\operatorname{Hom}}
\newcommand{\io}{\iota}
\newcommand{\oY}{\overline{Y}}
\newcommand{\ox}{\overline{x}}

\newcommand{\eps}{\epsilon}
\newcommand{\dimq}{\operatorname{dim}_q}
\newcommand{\ttil}{\tilde{t}}
\newcommand{\om}{\omega}
\newcommand{\Om}{\Omega}
\newcommand{\Omtil}{\widetilde{\Omega}}
\newcommand{\Tr}{\operatorname{Tr}}
\newcommand{\multq}{\operatorname{mult}_q}
\newcommand{\mult}{\operatorname{mult}}

\newcommand{\R}{\mathbb{R}}
\newcommand{\C}{\mathbb{C}}
\newcommand{\la}{\langle}
\newcommand{\ra}{\rangle}
\newcommand{\cB}{\mathcal{B}}
\newcommand{\cA}{\mathcal{A}}
\newcommand{\cC}{\mathcal{C}}
\newcommand{\cst}{C$^*$}

\newcommand{\Ah}{\hat{A}}
\newcommand{\deh}{\hat{\Delta}}
\newcommand{\N}{\mathbb{N}}
\newcommand{\NstarN}{\N \star \N}
\newcommand{\al}{\alpha}
\newcommand{\be}{\beta}
\newcommand{\GL}{\operatorname{GL}}

\newcommand{\cK}{\mathcal{K}}
\newcommand{\meq}{\underset{\text{\rm mon}}{\sim}}
\newcommand{\SU}{\operatorname{SU}}
\newcommand{\cstr}{C^*_r}
\newcommand{\cstu}{C^*_u}
\newcommand{\Vh}{\hat{V}}
\newcommand{\oF}{\overline{F}}
\newcommand{\length}{\operatorname{length}}
\newcommand{\tot}{\mathbin{\text{\footnotesize\textcircled{\tiny \sf T}}}}
\newcommand{\M}{\operatorname{M}}
\newcommand{\U}{\operatorname{U}}
\newcommand{\cF}{\mathcal{F}}
\newcommand{\cJ}{\mathcal{J}}

\begin{document}
\begin{center}
{\LARGE\bf Ergodic coactions with large multiplicity \\ and monoidal
equivalence of quantum groups}

\bigskip

{\sc by Julien Bichon$^{\text{\rm (a)}}$, An De Rijdt$^{\text{\rm (b)}}$ and Stefaan Vaes$^{\text{\rm (b,c)}}$}
\end{center}

{\footnotesize (a) \parbox[t]{15cm}{Laboratoire de Math{\'e}matiques Appliqu{\'e}es; Universit{\'e} de Pau et des Pays de l'Adour; IPRA;\\ Avenue de
l'Universit{\'e}, F--64000 Pau
(France)} \\
(b) Department of Mathematics; K.U.Leuven; Celestijnenlaan 200B; B--3001 Leuven (Belgium) \\
(c) Institut de Math{\'e}matiques de Jussieu; Alg{\`e}bres d'Op{\'e}rateurs; 175, rue du Chevaleret; F--75013 Paris (France)

e-mail: julien.bichon@univ-pau.fr, an.derijdt@wis.kuleuven.ac.be, vaes@math.jussieu.fr}

\begin{abstract}
\noindent We construct new examples of ergodic coactions of compact
quantum groups, in which the multiplicity of an irreducible
corepresentation can be strictly larger than the dimension of the
latter. These examples are obtained using a bijective correspondence
between certain ergodic coactions on \cst-algebras and unitary fiber functors on the
representation category of a compact quantum group. We classify these unitary fiber functors on
the universal orthogonal and unitary quantum groups. The associated
\cst-algebras and von Neumann algebras can be defined by generators
and relations, but are not yet well understood.
\end{abstract}

\section*{Introduction}

By a well known theorem of H\o egh-Krohn, Landstad  and St\o rmer \cite{HKLS}, compact groups only admit ergodic actions on tracial \cst-algebras.
Indeed, the (unique) invariant state is necessarily a trace. Moreover, given a compact group $G$ acting ergodically on a \cst-algebra $B$, one
studies the so-called \emph{spectral subspaces}: the action of $G$ on $B$ yields a unitary representation of $G$ which can be decomposed into a
direct sum of irreducible representations. One proves that the multiplicity of an irreducible representation is necessarily bounded by the dimension
of this irreducible representation.

A deeper analysis of the spectral structure of ergodic actions of compact groups has been made by A.\ Wassermann \cite{wasser1,wasser2,wasser3}. In
the culmination of his work, Wassermann shows that the compact group $\SU(2)$ only admits ergodic actions on von Neumann algebras of finite type I.

In the 1980's, Woronowicz introduced the notion of a \emph{compact quantum group} and generalized the classical Peter-Weyl representation theory.
Many fascinating examples of compact quantum groups are available by now: Drinfel'd and Jimbo \cite{drinfeld,jimbo} introduced the
\emph{$q$-deformations of compact semi-simple Lie groups} and Rosso \cite{rosso} showed that they fit into the theory of Woronowicz. The
\emph{universal orthogonal and unitary} quantum groups were introduced by Van Daele and Wang \cite{VDW} and studied in detail by Banica
\cite{banica2,banica1}. Other examples of compact quantum groups, related with graphs and metric spaces have been constructed by the first author
\cite{bichon4} and by Banica \cite{banica3,banica4}.

The abstract theory of \emph{ergodic coactions} of compact quantum groups on \cst-algebras has been initiated by Boca \cite{boca} and Landstad
\cite{landstad}. The major difference with the compact group case, is the following: the (unique) invariant state is no longer a trace. Indeed, Wang
\cite{wang} gave examples of ergodic coactions of universal unitary quantum groups on type III factors.

Moreover, in the work of Boca, the multiplicity of an irreducible corepresentation is proved to be bounded by the \emph{quantum dimension} rather
than the ordinary dimension of the corepresentation. Nevertheless, in
all the examples known up to now, the multiplicity is actually bounded by the ordinary
dimension. In this paper, we provide examples of ergodic coactions where the multiplicity of an irreducible corepresentation is strictly larger than
the ordinary dimension of the corepresentation.

In \cite{wasser3}, A.\ Wassermann gives a complete classification of the ergodic actions of $\SU(2)$, essentially labeling them by the finite
subgroups of $\SU(2)$. It would, of course, be great to give a complete classification of ergodic coactions of the deformed $\SU_q(2)$. In
\cite{tomatsu}, Tomatsu provides a first step in this direction: he computes all ergodic coactions of $\SU_q(2)$ on \lq virtual\rq\ quotient spaces
$\SU_q(2)/\Gamma$. (More precisely, he describes all the coideals of the quantum group $\SU_q(2)$.) By construction, the ergodic coactions of
$\SU_q(2)$ on its virtual quotient spaces are such that the multiplicity of an irreducible corepresentation is bounded by its dimension. The results
of this paper imply in particular that there are much more ergodic coactions of $\SU_q(2)$ than the ones studied by Tomatsu.

The major tool to produce our new examples of ergodic coactions of compact quantum groups, is the notion of \emph{monoidal equivalence} of quantum
groups. One can look at a compact quantum group with several degrees of precision. At first, we study only the \emph{fusion rules} in the
representation theory: we label the irreducible corepresentations and describe how a tensor product of irreducibles breaks up into irreducibles.
Taking into account only these fusion rules, we loose a lot of information: for example, the $q$-deformed compact Lie groups have the same fusion
rules as their classical counterparts. In a next approximation, one studies the corepresentation theory of a compact quantum group as a
\emph{monoidal category}, but without its concrete realization (the so-called \emph{forgetful functor} to the category of Hilbert spaces). This is
crucial: by the Tannaka-Krein reconstruction theorem \cite{wor3}, the concrete monoidal category of (finite-dimensional) corepresentations
essentially determines the compact quantum group. Note that knowing the representation theory of a compact quantum group as a monoidal category,
comes down to knowing the \emph{fusion rules together with the $6j$-symbols}, see Remark \ref{rem.6j}.

Closely related to the notion of monoidally equivalent quantum groups, is the notion of a \emph{unitary fiber functor} on a compact quantum group.
Essentially, a unitary fiber functor gives \emph{another concrete realization}, different from the tautological realization, of the representation
theory of a compact quantum group.

In this paper, we choose not to use the abstract language of categories. We give \lq down-to-earth definitions\rq\ of monoidally equivalent quantum
groups and unitary fiber functors, see \ref{def.moneq} and \ref{def.fiber}. This makes the construction of associated \cst-algebras and coactions
straightforward. This concrete approach is well adapted to the language of corepresentations of compact quantum groups. In this way, using previous
results of Banica \cite{banica2,banica1}, we can show very easily as well the monoidal equivalence of the universal orthogonal and universal unitary
quantum groups.

The results in this paper can be summarized as follows.

\begin{list}{$\bullet$}{\setlength{\labelwidth}{3mm}\setlength{\leftmargin}{4mm}}
\item In Section \ref{sec.spectral}, we recall the theory of spectral subspaces \cite{boca,landstad} and provide a simple proof for the multiplicity
bound. We also introduce the notion of \emph{quantum multiplicity} of
an irreducible corepresentation in an ergodic coaction and this can be strictly
larger than the ordinary multiplicity.
\item In Theorem \ref{thm.main}, we show that there is a natural bijective correspondence between certain ergodic coactions of compact quantum groups
and unitary fiber functors. These coactions are called \emph{of full quantum multiplicity}. These are precisely the ergodic coactions for which the
crossed product is isomorphic with the compact operators, see \cite{landstad}. They can also be described as the ergodic coactions for which the
inequality between quantum multiplicity and quantum dimension, becomes an equality.
\item In Section \ref{sec.cocycle}, we study the special case of unitary fiber functors preserving the dimension. This leads to a bijective
correspondence with \emph{unitary $2$-cocycles on the dual, discrete, quantum group}. The ideas for this section come from the work of Wassermann
\cite{wasser2}.
\item In Sections \ref{sec.ao}, we establish the monoidal equivalence between the universal orthogonal quantum groups $A_o(F)$. Recall that,
for any $F \in \GL(n,\C)$ satisfying $F \overline{F} = \pm 1$, one defines the compact
quantum group $A_o(F)$ as the universal quantum group generated by the coefficients of a unitary $n$ by $n$ matrix $U$ with relations $U = F
\overline{U} F^{-1}$. The comultiplication on $A_o(F)$ is (uniquely)
defined in such a way that $U$ becomes a corepresentation.
We show that $A_o(F_1)$ is monoidally equivalent with $A_o(F_2)$ if
and only if the signs of the $F_i \overline{F_i}$ agree and
$\Tr(F_1^* F_1) = \Tr(F_2^* F_2)$. In particular, if $0 < q \leq 2 - \sqrt{3}$, there is a continuous family of non-isomorphic $A_o(F)$ monoidally
equivalent with $\SU_q(2)$.
\item  In Section \ref{sec.au}, we prove similar results for the universal unitary quantum groups $A_u(F)$, defined as the universal quantum group
generated by the coefficients of a unitary $n$ by $n$ matrix $U$ with
the relation that $F \overline{U} F^{-1}$ is unitary. Again, the
comultiplication is defined such that $U$ becomes a corepresentation.
We show that the quantum dimension of $U$, i.e.\
$\sqrt{\Tr(F^*F) \Tr((F^*F)^{-1})}$, is a complete invariant for the $A_u(F)$ up to monoidal equivalence.
\item Using the previous results, we obtain a complete classification of the ergodic coactions of full
quantum multiplicity of $A_o(F)$ and $A_u(F)$, as well as a computation of the $2$-cohomology of their duals (Corollary \ref{cor.classao}). In
particular, we construct ergodic coactions of $\SU_q(2)$ such that the multiplicity of the fundamental corepresentations is arbitrarily large
(Corollary \ref{cor.largemult}).
\end{list}

In the theory of Hopf algebras, ergodic coactions of full quantum multiplicity correspond to \emph{Hopf-Galois extensions}. In this algebraic
setting, several results related to ours have been obtained. The relation between Hopf-Galois extensions and fiber functors is due to Ulbrich
\cite{ulbrich} and the relation between monoidal equivalence of Hopf algebras and Hopf-bi-Galois extensions has been established by Schauenburg
\cite{schauenburg}. Fiber functors preserving the dimension and $2$-cocycles have been studied by Etingof and Gelaki \cite{etingof-gelaki}. The main
difference between these Hopf algebraic results and our work, lies in dealing with the $^*$-structure and positivity. In a sense, we are dealing with
the real forms of (certain) Hopf algebras. This allows us to construct Hilbert space representations and C$^*$-algebras. The compatibility of fiber
functors with $^*$-structures is a severe restriction. Indeed, there exist many fiber functors on the representation category of $\SU(2)$ (see
\cite{bruguieres}), but the forgetful functor is the only one compatible with the $^*$-structure.

After completion of a first version of this paper -- signed by the last two authors -- the first author joined the project and his preprint
\cite{bichon1} was taken into account, yielding the current paper as a final result.

\section{Preliminaries}
\begin{notation}
Consider a subset $S$ of a \cst-algebra. We denote by $\langle S \rangle$
the linear span of $S$ and by $[S]$ the closed linear span of $S$.

The symbol $\ot$ denotes \emph{tensor products} of Hilbert spaces, \emph{minimal}
tensor products of \cst-algebras as well as algebraic tensor products
of $^*$-algebras.

If $A$ is a $^*$-algebra and $U \in M_{n_1,n_2}(\C) \ot A$, we denote by $\overline{U}$ the matrix $\overline{U}_{ij} = U_{ij}^*$.

We make use of the \emph{leg numbering notation}. For instance, if $v \in A \ot B$,
then $v_{12}$ denotes the element in $A \ot B \ot C$ defined by
$v_{12} = v \ot 1$. We analogously use the notations $v_{13}$,
$v_{23}$, etc.
\end{notation}

\begin{definition}
A \emph{compact quantum group} $\cG=(A,\de)$ consists of
a unital \cst-algebra $A$ together with a unital *-homomorphism
$\de:A\to A\ot A$ satisfying the coassociativity relation
\[(\io\ot\de)\de=(\de\ot\io)\de\] and the cancellation properties
\[[\de(A)(A\ot 1)]=A\ot A=[\de(A)(1\ot A)]\;.\]
\end{definition}

If $\cG=(A,\de)$ is a compact quantum group, there exists a unique
state $h$ on $A$ which is invariant under the comultiplication:
\[(\io\ot h)\de(a)=(h\ot\io)\de(a)=h(a)1\] for all $a\in A$. We call
$h$ the \emph{Haar state} of $\cG$.

\begin{definition}
Let $H$ be a Hilbert space. A \emph{unitary corepresentation} $v$ of
$\cG$ on $H$ is a unitary element of $\M(\cK(H)\ot A)$ satisfying
$(\io\ot\de)(v)=v_{12}v_{13}$. The dimension of the underlying Hilbert space $H$ is called the dimension of $v$ and denoted by $\dim v$.
\end{definition}

The tensor product of the unitary corepresentations $v$ and $w$ is
defined by
$$v \tot w := v_{13} w_{23} \; .$$
Note that the corepresentations $v \tot w$ and $w \tot v$ are in general not unitarily equivalent. This is a crucial feature of quantum groups.

\begin{notation}
Let $\cG=(A,\de)$ be a compact quantum group. Given two unitary
corepresentations $v\in \M(\cK(H)\ot A)$ and $w\in \M(\cK(K)\ot A)$, we
denote by $\Mor(v,w)$ the intertwiners between $v$ and $w$:
\[\Mor(v,w)=\{S\in \B(H,K)\mid (S\ot 1)v=w(S\ot 1)\}\;.\]
\end{notation}

\begin{terminology}
A unitary corepresentation $v\in \M(\cK(H)\ot A)$ is called \emph{irreducible}
if $\Mor(v,v)=\C 1$. A unitary corepresentation $w\in \M(\cK(L)\ot A)$ is called
\emph{unitarily equivalent} to $v$ if $\Mor(v,w)$ contains a unitary operator.
\end{terminology}

In this paper, all corepresentations are assumed to be unitary.
Recall the following well known facts (see \cite{wor2}).
Every irreducible corepresentation of a compact quantum group is
finite dimensional and every corepresentation decomposes as a direct
sum of irreducible corepresentations.

\begin{notation}
Let $\cG=(A,\de)$ be a compact quantum group. We denote by $\cGh$ the
set of equivalence classes of irreducible corepresentations of $\cG$
and we choose
unitary representatives $U^x \in \B(H_x) \ot A$ for all $x \in
\cGh$. We denote by $\eps$ the (class of the) trivial corepresentation
$1 \in A$.
\end{notation}

For every $x\in \cGh$, there is a unique $\ox\in\cGh$ such that $\Mor(\eps,x\tot\ox) \neq 0$. The irreducible corepresentation $\ox$ is called the
\emph{adjoint} of $x$. Take now $t\in\Mor(\eps,x\tot\ox)$, $t\neq 0$ and define the antilinear map \[S_t:H_x\to H_{\ox}:\xi\mapsto (\xi^*\ot 1)t\;.\]
Define $Q_x:=S_t^*S_t$ and normalize $t$ in such a way that $\Tr(Q_x)=\Tr(Q_x^{-1})$. This uniquely determines $Q_x$ and fixes $t$ up to a number of
modulus $1$. If we take the unique $\ttil\in\Mor(\eps,\ox\tot x)$ such that $(t^*\ot 1)(1\ot \ttil)=1$, then $S_{\ttil}=S_t^{-1}$ and
$Q_{\ox}=(S_tS_t^*)^{-1}$.

\begin{notation}
The \emph{quantum dimension} of $x \in \cGh$ is defined as $\Tr(Q_x)$
and denoted by $\dimq(x)$.
\end{notation}

Observe that $t^*t=\dimq(x)1$ and that $\dimq(\ox)=\dimq(x)\geq\dim(x)$. The orthogonality relations can then be written as follows (see
\cite{wor2}): for $\xi\in H_x$, $\eta\in H_y$,
\begin{equation}\label{eq.haarstate}
(\io\ot h)(U^x(\xi\eta^*\ot 1)(U^y)^*)=\frac{\sde_{x,y} 1}{\dimq(x)}\langle Q_x\xi,\eta\rangle\;.
\end{equation}

\begin{notation}
Let $\cG=(A,\de)$ be a compact quantum group. We
denote by $\cA$ the set of coefficients of finite dimensional
corepresentations of $\cG$. Hence,
\[\cA =\la (\om_{\xi,\eta}\ot\io)(U^x)\mid x\in\cGh,\ \xi,\eta\in H_x
\ra \]
Then, $\cA$ is a unital dense $^*$-subalgebra of $A$. Restricting $\de$ to
$\cA$, $\cA$ becomes a Hopf $^*$-algebra.
\end{notation}

\begin{terminology}
Let $\cB$ be a unital $^*$-algebra. A linear functional $\om : \cB
\recht \C$ is said to be a \emph{faithful state} on $\cB$ if $\om(1) = 1$ and
if $\om(aa^*) \geq 0$ for all $a \in \cB$, with equality holding if
and only if $a=0$.
\end{terminology}

Observe that the Haar state of a compact quantum group $(A,\de)$ is a faithful state on the underlying Hopf $^*$-algebra $\cA$.

The dual (discrete) quantum group of $(A,\de)$, can be defined as
$$\Ah = \bigoplus_{x \in \cGh} \B(H_x) \quad\text{with}\quad \deh(a) S = S a \quad\text{for all}\;\; a \in \Ah, S \in \Mor(x,y \tot z) \; .$$

There is, of course, another way to define $(\Ah,\deh)$, using the
left regular corepresentation of $(A,\de)$ and the theory of
multiplicative unitaries. We only need this in Section
\ref{sec.cocycle}, see Notation \ref{not.multun}.

\section{Spectral subspaces and quantum multiplicity}
\label{sec.spectral}

In this section, we give a brief overview of the general theory of ergodic coactions of compact quantum groups. We study spectral subspaces and prove
in particular that they are finite dimensional. The results in this section are well known (see \cite{HKLS} for the classical case of compact groups
and \cite{boca,landstad,tomatsu} for compact quantum groups). We give a short presentation for the convenience of the reader.

Let $\cG = (A,\de)$ be a compact quantum group. Recall that $\cGh$ denotes the set of equivalence classes of irreducible corepresentations of $\cG$
and that we chose unitary representatives $U^x \in \B(H_x) \ot A$ for all $x \in \cGh$. We have $(\io \ot \de)(U^x) = U^x_{12} U^x_{13}$.

\begin{definition} Let $B$ be a unital C*-algebra. A (right) coaction
  of $(A,\de)$ on $B$ is a unital $^*$-homomorphism
$\sde : B \recht B \ot A$ satisfying
\[(\sde\ot \io)\sde = (\io \ot \de)\sde \quad\text{and}\quad [\sde(B)(1\ot A)]=B\ot A\; .\] The
coaction $\sde$ is said to be ergodic if the fixed point algebra $B^{\sde}:=\{x\in
A\mid  \sde (x)=x\ot 1\}$ equals $\C 1$.
\end{definition}

\begin{remark}
If $\sde : B \recht B \ot A$ is an ergodic coaction of $(A,\de)$ on
$B$ there is a unique state $\om$ on $B$ which is invariant under
$\sde$, given by $\om(b) 1 = (\io
\ot h)\sde(b)$.
\end{remark}

In what follows, we fix a coaction $\sde:B\to B\ot A$ of a compact
quantum group $\cG=(A,\de)$.

\begin{definition}
Let $x \in \cGh$. We define the spectral
subspace associated with $x$ by
$$K_x = \{X \in \oH_x \ot B \mid (\io \ot \sde)(X) = X_{12} U^x_{13} \} \;
.$$
\end{definition}

Defining $\Hom(H_x,B)= \{
S:H_x \to B\mid S \ \textrm{linear and}\ \sde(S \xi )=(S \ot \io
)(U^x(\xi \ot 1))\}$, we have $K_x \cong \Hom(H_x,B)$, associating
to every
$X\in K_x$ the operator $S_X:H_x\to B:\xi\mapsto X(\xi\ot 1)$.

\begin{definition} \label{def.densesub}
We define $\cB$ as the subspace of $B$ generated by the spectral
subspaces, i.e.\
$$\cB:= \la X(\xi\ot 1)\mid x\in\cGh,\ X\in K_x,\ \xi\in
H_x \ra \; .$$
\end{definition}

Observe that $\cB$ is a dense unital *-subalgebra of $B$ and that the
restriction
$\sde:\cB\to\cB\ot\cA$ defines a coaction of the Hopf $^*$-algebra
$(\cA,\de)$ on $\cB$.

\begin{terminology}
A coaction $\sde:B\to B\ot A$ of $(A,\de)$ on $B$ is said to be
\emph{universal} if $B$ is the universal
enveloping \cst-algebra of $\cB$. It is said to be \emph{reduced}
if the conditional expectation $(\io\ot h)\sde$ of $B$ on $B^\sde$ is faithful.
\end{terminology}

\begin{remark} \label{rem.coamenable}
Observe that an ergodic coaction is reduced if and only if the unique invariant state is faithful.

In the special case where $B = A$ and $\sde = \de$, the $^*$-algebra $\cB$ coincides with the underlying Hopf $^*$-algebra $\cA \subset A$ consisting
of coefficients of finite-dimensional corepresentations. So, we obtain the usual notions: a compact quantum group $(A,\de)$ is said to be universal
if $A$ is the universal enveloping \cst-algebra of $\cA$ and reduced if the Haar state is faithful on $A$. Clearly, any compact quantum group has its
universal and reduced companion. In the case where $(A,\de)$ is the dual of a discrete group, these notions coincide with the full, resp.\ reduced
group \cst-algebra.

A compact quantum group is said to be \emph{co-amenable} its universal and reduced companion coincide. Equi\-valently, a compact quantum group is
co-amenable if we have on the same \cst-algebra a bounded co-unit and a faithful Haar state. It is then clear that a coaction of a co-amenable
compact quantum group is always both universal and reduced. Examples of co-amenable compact quantum groups include $\SU_q(n)$ and other
$q$-deformations of compact Lie groups.
\end{remark}

Fix now an \emph{ergodic} coaction $\sde : B \recht B \ot A$ of a
compact quantum group $(A,\de)$ on a unital \cst-algebra $B$. Denote
by $\om$ the unique invariant state on $B$.

For every $X,Y\in
K_x$, $XY^*\in B^{\sde}=\C 1$. So, for $x \in \cGh$, the spectral subspace $K_x$
is a Hilbert space with scalar product $\la X,Y \ra 1 = XY^*$.

\begin{terminology}
Let $\sde$ be an ergodic coaction of $(A,\de)$ on $B$. The dimension
of the Hilbert space $K_x$ is called the \emph{multiplicity of $x$ in
  $\sde$} and denoted by $\mult(x)$.
\end{terminology}

Define, for every $x \in \cGh$, the element $X^x \in \B(H_x,\oK_x)
\ot B$ such that $(X^x)^*(\oY \ot 1) = Y^*$ for all $Y \in K_x$.
Observe that $X^x (X^x)^* = 1$. Therefore, $(X^x)^* X^x\in \B(H_x)\ot
B$ is a projection.

Let $x \in \cGh$. Take $t \in \Mor(\eps,\ox \tot x)$, normalized in such a way that $t^* t = \dimq(x)$. Define the antilinear map
$$T_t : K_x \recht K_{\ox} : T_t(Y) = (t^* \ot 1)(1 \ot Y^*) \; .$$
Since $t$ is fixed up to a number of modulus one, $L_x := T_t^* T_t$
is a well defined positive element of $\B(K_x)$.

\begin{definition}
We put $\multq(x) := \sqrt{\Tr(L_x) \Tr(L_{\ox})}$ and call $\multq(x)$ the \emph{quantum
  multiplicity of $x$ in $\sde$}.
We prove in Theorem \ref{thm.fullquantum} that $\multq(x) \leq
\dimq(x)$ for all $x \in \cGh$. If equality holds for all $x \in
\cGh$, we say that $\sde$ is of \emph{full quantum multiplicity}.
\end{definition}

\begin{theorem} \label{thm.fullquantum}
Let $\sde : B\to B\ot A$ be an ergodic coaction of a compact quantum
group $\cG=(A,\de)$ on a unital \cst-algebra $B$.
For every irreducible corepresentation $x\in\cGh$, \[\mult(x)
\leq \multq(x) \leq
  \dimq(x) \;.\] With $X^x$ defined as above, the ergodic coaction
  is of full quantum multiplicity if and only if $X^x$ is unitary
  for all $x\in\cGh$.
\end{theorem}

\begin{proof}
Let $x \in \cGh$ and take $t \in \Mor(\eps,\ox \tot x)$, normalized in such a way that $t^* t = \dimq(x)$. If we take the unique $\ttil \in
\Mor(\eps,x \tot \ox)$ satisfying $(1 \ot t^*)(\ttil \ot 1) = 1$, we have $\ttil^* \ttil = \dimq(\ox)$ and hence, it is clear that $T_{\ttil} =
T_t^{-1}$ and $L_{\ox} = (T_t T_t^*)^{-1}$. So,
$$\multq(x) = \sqrt{\Tr(T_t T_t^*) \Tr((T_t T_t^*)^{-1})} \geq
\dim(T_t) = \mult(x) \; ,$$
for all $x \in \cGh$. By definition
\begin{equation}\label{eq.hulphulp}
\la L_x X, Y \ra 1 = \la (t^* \ot 1)(1 \ot Y^*) , (t^* \ot 1)(1 \ot
X^*) \ra 1 = (t^* \ot 1)(1 \ot Y^*X) (t \ot 1) \; .
\end{equation}
Since $t^*(1 \ot a)t = \Tr(Q_x^{-1} a)$ for all $a \in \B(H_x)$, we
conclude that
$$\Tr(L_x)1 = (\Tr(Q_x^{-1} \cdot) \ot \io)((X^x)^* X^x) \; .$$
Since $(X^x)^* X^*$ is a projection, it follows that $\Tr(L_x) \leq
\Tr(Q_x^{-1}) = \dimq(x)$ for all $x \in \cGh$. Applying this
inequality to $x$ and $\ox$, we conclude that $\multq(x) \leq
  \dimq(x)$ for all $x \in \cGh$.

Moreover, $\multq(x)=\dimq(x)$ for all $x \in \cGh$ if and only if $(\Tr(Q_x^{-1} \cdot) \ot \io)((X^x)^* X^x) = \Tr(Q_x^{-1})1$ for all $x \in
\cGh$. This last statement holds if and only if $(X^x)^* X^x = 1$ for all $x \in \cGh$, i.e.\ when $X^x$ is unitary for all $x \in \cGh$.
\end{proof}

It is also straightforward to show that the invariant state $\om$ is a KMS state, at least for universal and for reduced ergodic coactions.

\begin{proposition}
If the ergodic coaction $\sde$ of $(A,\de)$ on $B$ is either universal
or reduced, the invariant state $\om$ is a KMS state. The elements of
the dense $^*$-subalgebra $\cB \subset B$ (Definition \ref{def.densesub}) are analytic with respect to
the modular group, given by
$$\sigma_t^\om \bigl( Y(\xi \ot 1) \bigr) = (L_x^{-it} Y) ( Q_x^{it}
\xi \ot 1) \quad\text{for all}\;\; x \in \cGh, Y \in K_x, \xi \in H_x
\; .$$
\end{proposition}

\begin{proof}
Let $y,z \in \cGh$ and $Y \in K_y, Z \in K_z$. One verifies that $(U^z)^*((\io \ot \om)(Z^* Y)\ot 1)U^y= (\io \ot \om)(Z^* Y)\ot 1$. Hence, $(\io \ot
\om)(Z^*Y)$ equals $0$ if $y \neq z$ and is scalar if $y=z$. Applying $\om$ to \eqref{eq.hulphulp}, we conclude that
$$(\io \ot \om)(Z^* Y) = \frac{\sde_{y,z} 1}{\dimq(y)} \la L_y Y, Z \ra \; .$$
If moreover $\xi \in H_y, \eta \in H_z$ and if we put $a = Y(\xi \ot 1)$ and $b = Z(\eta \ot 1)$, we get
\begin{equation}\label{eq.invariant-state}
\om(b^*a)=\frac{\sde_{y,z}}{\dimq(y)}\langle\xi,\eta\rangle\langle L_y Y,Z\rangle\;.
\end{equation}
Using \eqref{eq.haarstate}, one checks that
$$\om(ab^*) = (\io \ot \om)(Y(\xi \eta^* \ot 1)Z^*) = \frac{\sde_{y,z} 1}{\dimq(y)} \la Q_y \xi,\eta \ra \; \la Y,Z \ra \; .$$
As a linear space, $\cB \cong \bigoplus_{x \in \cGh} (K_x \ot H_x)$. So, we can define linear maps $\sigma_t^\om : \cB \recht \cB$ by the formula
$$\sigma_t^{\om}\bigl(Y (\xi \ot 1)
\bigr):=L_y^{-it}(Y)(Q_y^{it}\xi \ot 1) \; .$$ It is clear that $(\sigma_t^\om)$ is a one-parameter group of linear isomorphisms of $\cB$. Observe
that all elements of $\cB$ are analytic with respect to $(\sigma_t^\om)$ and that $\om(\sigma_i^\om(a)b^*) = \om(b^*a)=\om(a\sigma^{\om}_{i}(b)^*)$
for all $a,b \in \cB$. Since \eqref{eq.invariant-state} implies that $\om$ is faithful on $\cB$, it follows that $\sigma^\om_{i} : \cB \recht \cB$ is
multiplicative and that $\sigma^\om_i(a)^* = \sigma^\om_{-i}(a^*)$ for all $a \in \cB$. Standard complex analysis allows to conclude that the
$\sigma^\om_t$ are $^*$-automorphisms of $\cB$. It is also clear that $\om$ is invariant under $\sigma_t^\om$.

If $\sde$ is a universal coaction, the one-parameter group $(\sigma_t^\om)$ extends to $B$ by universality. If $\sde$ is a reduced coaction, we can
extend $\sigma_t^\om$ to $B$ because $\om$ is invariant under $\sigma_t^\om$ and $\om$ is faithful on $B$. In both cases, it follows that
$(\sigma_t^\om)$ satisfies the KMS condition with respect to $\om$ and so, $\om$ is a KMS state.
\end{proof}

Finally observe that $\om$ is a trace if and only if
$L_x=1$ and $Q_x=1$ for all $x\in\cG$ with $K_x \neq 0$.

\section{Commuting coactions and monoidal equivalence of quantum
  groups}

Our main goal in this section is to show the relation
between ergodic coactions of full quantum multiplicity and
monoidal equivalence of compact quantum groups. In Sections
\ref{sec.ao} and \ref{sec.au}, we shall give examples of monoidally
equivalent compact quantum groups, giving rise to new examples of
ergodic coactions.

The relation between ergodic coactions of full quantum multiplicity
and monoidal equivalence allows us to classify completely such
coactions for the unitary and orthogonal quantum groups $A_u(F)$ and
$A_o(F)$, in particular for $\SU_q(2)$.

\begin{definition} \label{def.moneq}
Two compact quantum groups $\cG=(A,\de)$ and $\cG_2=(A_2,\de_2)$ are said to be \emph{monoidally equivalent} if there exists a bijection
$\vphi:\cGh\to\cGh_2$ satisfying $\vphi(\eps) = \eps$, together with linear isomorphisms
\[\vphi:\Mor(x_1 \tot \cdots \tot x_r ,y_1\tot\cdots\tot
y_k)\to\Mor(\vphi(x_1) \tot \cdots \tot \vphi(x_r),\vphi(y_1)\tot\cdots \tot \vphi(y_k))\] satisfying the following conditions:
\begin{equation}\label{eq.moneq}
\begin{alignedat}{2}
\vphi(1) &= 1 & \qquad
\vphi(S \ot T) &= \vphi(S) \ot \vphi(T) \\
\vphi(S^*) &= \vphi(S)^* & \qquad
\vphi(S T) &=\vphi(S) \vphi(T)
\end{alignedat}
\end{equation}
whenever the formulas make sense. In the first formula, we consider $1 \in \Mor(x,x) = \Mor(x, x \tot \eps) = \Mor(x,\eps \tot x)$. Such a collection
of maps $\vphi$ is called a \emph{monoidal equivalence} between $\cG$ and $\cG_2$.
\end{definition}

\begin{remark} \label{rem.reduction}
To define a monoidal equivalence between $\cG$ and $\cG_2$, it suffices to define a bijection $\vphi:\cGh\to\cGh_2$ satisfying $\vphi(\eps) = \eps$,
together with linear isomorphisms $\vphi : \Mor(x,y_1 \tot \cdots \tot y_k) \recht \Mor(\vphi(x),\vphi(y_1) \tot \cdots \tot \vphi(y_k))$ for
$k=1,2,3$, satisfying
\begin{alignat}{2}
\vphi(1) &= 1 & & \label{eq.unital} \\
\vphi(S)^* \vphi(T) &= \vphi(S^*T) & &\quad\text{for all}\quad S \in
\Mor(a,x \tot y), T \in \Mor(b,x \tot y) \label{eq.fusion} \\
\vphi((S \ot 1)T) &= (\vphi(S) \ot 1)\vphi(T) & &\quad\text{for
  all}\quad T \in \Mor(a,b \tot z), S \in \Mor(b,x \tot y)
\label{eq.6j1} \\
\vphi((1 \ot S)T) &= (1 \ot \vphi(S)) \vphi(T) & &\quad\text{for
  all}\quad T \in  \Mor(a,x \tot b), S \in \Mor(b,y \tot z) \label{eq.6j2}
\end{alignat}
Indeed, such a $\vphi$ admits a unique extension to a monoidal equivalence. Again, \eqref{eq.unital} should be valid for $1 \in \Mor(x,x) = \Mor(x, x
\tot \eps) = \Mor(x,\eps \tot x)$.
\end{remark}

\begin{remark} \label{rem.6j}
Observe that the existence of the linear isomorphisms $\vphi : \Mor(a,b \tot c) \recht \Mor(\vphi(a),\vphi(b) \tot \vphi(c))$ only says that $\cG$
and $\cG_2$ have the \emph{same fusion rules}. Adding \eqref{eq.unital}--\eqref{eq.6j2} means that $\cG$ and $\cG_2$ moreover have the \emph{same
  $6j$-symbols} (see \cite{CFS}).
Indeed, taking orthonormal bases for all $\Mor(a,b \tot c)$, we can write two natural orthonormal bases for $\Mor(a,x \tot y \tot z)$, one given by
elements $(S \ot 1)T$, the other given by elements $(1 \ot S)T$. The coefficients of the transition unitary between both orthonormal bases are called
the $6j$-symbols of $\cG$.
\end{remark}

\begin{remark}
As we shall see in Sections \ref{sec.ao} and \ref{sec.au}, there are
natural examples of monoidal equivalences where $\dim(\vphi(x)) \neq
\dim(x)$. On the other hand, it is clear that $\dimq(\vphi(x)) =
\dimq(x)$ for all $x \in \cG$ and all monoidal equivalences $\vphi$.
We shall see in Section \ref{sec.au} that for a certain class of
compact quantum groups (the universal unitary ones), this equality of
quantum dimension is the only constraint for monoidal equivalence.
\end{remark}

\begin{remark}
It is clear that a monoidal equivalence $\vphi$ in the sense of Definition \ref{def.moneq}, defines a monoidal equivalence in the usual sense
(preserving the $^*$-operation), between the monoidal categories of finite dimensional corepresentations of $\cG$ and $\cG_2$ (\cite{maclane}).
Moreover, this monoidal equivalence is uniquely determined up to isomorphism. We prefer to work with the \lq concrete\rq\ data of Definition
\ref{def.moneq}, avoiding all kinds of identifications.
\end{remark}

\begin{notation}
If two compact quantum groups $\cG=(A,\de)$ and $\cG_2=(A_2,\de_2)$
are monoidally equivalent, we write $\cG\meq\cG_2$
\end{notation}

Closely related to the notion of monoidal equivalence, is the
following notion of \emph{unitary fiber functor} (see Proposition
\ref{prop.tannaka-krein} for the relation between both notions).

\begin{definition} \label{def.fiber}
Let $\cG=(A,\de)$ be a compact quantum group. A \emph{unitary fiber functor} associates to every $x \in \cGh$ a finite dimensional Hilbert space
$H_{\vphi(x)}$ and consists further of linear maps
$$\vphi : \Mor(x_1 \tot \cdots \tot x_r,y_1 \tot \cdots \tot y_k) \recht
\B(H_{\vphi(x_1)} \ot \cdots \ot H_{\vphi(x_r)},H_{\vphi(y_1)} \ot \cdots \ot H_{\vphi(y_k)})$$
satisfying equations \eqref{eq.moneq} in Definition \ref{def.moneq}.
\end{definition}

\begin{remark}
We make a remark analogous to \ref{rem.reduction}. To define a unitary fiber functor on $\cG$, it suffices to associate to every $x \in \cGh$ a
finite-dimensional Hilbert space $H_{\vphi(x)}$, with $H_{\vphi(\eps)} = \C$ and to define linear maps $\vphi : \Mor(x,y_1 \tot \cdots \tot y_k)
\recht \B(H_{\vphi(x)},H_{\vphi(y_1)} \ot \cdots \ot H_{\vphi(y_k)})$ for $k=1,2,3$, satisfying \eqref{eq.unital} -- \eqref{eq.6j2} as well as the
non-degenerateness assumption
$$\{\vphi(S) \xi \mid a \in \cGh, S \in \Mor(a,b \tot c), \xi \in
H_{\vphi(a)} \} \quad\text{is total in}\quad H_{\vphi(b)} \ot
H_{\vphi(c)}$$
for all $b,c \in \cGh$.

Moreover, it follows from Proposition \ref{prop.tannaka-krein} below, that a unitary fiber functor $\vphi$ on $\cG$ naturally defines a compact
quantum group $\cG_2$ such that $\vphi$ becomes a monoidal equivalence between $\cG$ and $\cG_2$.
\end{remark}

\begin{theorem} \label{thm.main}
Consider a compact quantum group $\cG=(A,\de)$ and let $\vphi$ be a unitary fiber functor on $\cG$.
\begin{itemize}
\item There exists a unique unital $^*$-algebra $\cB$ equipped with a faithful state
$\om$ and unitary elements $X^x \in \B(H_x,H_{\vphi(x)}) \ot \cB$ for
all $x \in \cGh$, satisfying
\begin{enumerate}
\item $X^y_{13} X^z_{23} (S \ot 1) = (\vphi(S) \ot 1)X^x
  \quad\text{for all}\quad S \in \Mor(x,y \tot z) \; ,$
\item the matrix coefficients of the $X^x$ form a linear basis of
  $\cB$,
\item $(\io \ot \om)(X^x) = 0 \quad\text{if}\quad x \neq \eps$.
\end{enumerate}
\item There exists a unique coaction $\sde : \cB \recht \cB \ot
\cA$ satisfying
$$(\io \ot \sde)(X^x) = X^x_{12} U^x_{13}$$
for all $x \in \cGh$.
\item The state $\om$ is invariant under $\sde$. Denoting
  by $B_r$ the \cst-algebra generated by $\cB$ in the
  GNS-representation associated with $\om$ and denoting by $B_u$ the
  universal enveloping \cst-algebra of $\cB$, the coaction $\sde$
  admits a unique extension to a coaction on $B_r$, resp.\ $B_u$.

These coactions are reduced, resp.\ universal and they are ergodic
and of full quantum multiplicity.
\item Every reduced, resp.\ universal, ergodic coaction of
  full quantum multiplicity, arises in this way from a unitary fiber functor.
\end{itemize}
\end{theorem}

\begin{proof}
Let $\vphi$ be a unitary fiber functor on $\cG=(A,\de)$. Define
the vector space $\cB = \oplus_{x \in \cGh} \B(H_x,H_{\vphi(x)})^*$.
We shall turn this vector space into a $^*$-algebra.

Define natural elements $X^x \in \B(H_x,H_{\vphi(x)}) \ot \cB$ by
$(\om_x\ot\io)(X^x)=(\sde_{x,y}\om_x)_{y\in\cGh}$ for all
$\om_x\in\B(H_x,H_{\vphi(x)})^*$. By definition, the coefficients of
the $X^x$ form a linear basis of $\cB$. Hence, it suffices to define a
product and an involution on the level of the $X^x$.

It is clear that there exists a unique bilinear multiplication map
$\cB \times \cB \recht \cB$ such that
$$X^y_{13} X^z_{23} (S \ot 1) = (\vphi(S) \ot 1)X^x
  \quad\text{for all}\quad S \in \Mor(x,y \tot z) \; .$$
But then
$$(X^a_{14} X^b_{24}) X^c_{34} \; ((S \ot 1) T \ot 1) = (\vphi(S) \ot
1 \ot 1) X^y_{13} X^c_{23} (T \ot 1) = ((\vphi(S) \ot 1)\vphi(T) \ot 1) X^x = (\vphi((S \ot 1)T) \ot 1) \; X^x$$ for all $S \in \Mor(y,a \tot b), T
\in \Mor(x,y \tot c)$. Since intertwiners of the form $(S\ot 1)T$ linearly span $\Mor(x,a \tot b \tot c)$, we conclude that
$$(X^a_{14} X^b_{24}) X^c_{34} \; (S \ot 1) = (\vphi(S) \ot 1) X^x$$
for all $S \in \Mor(x,a \tot b \tot c)$. Analogously,
$$X^a_{14} (X^b_{24} X^c_{34}) \; (S \ot 1) = (\vphi(S) \ot 1) X^x$$
for all $S \in \Mor(x,a \tot b \tot c)$. This proves the associativity of the product on $\cB$. It is clear that $X^\eps$ provides the unit element
of $\cB$.

Observe also that
\begin{equation}\label{fifu3}(\vphi(S)^*\ot 1)X^x_{13}X^y_{23}=X^z(S^*\ot
1)\end{equation} for all $S\in\Mor(z, x\tot y)$.

We define an antilinear map $b \mapsto b^*$ on $\cB$ such that
\begin{equation}\label{involution}
(X^x)^*_{13}(\vphi(t)\ot 1)=X^{\ox}_{23}(t\ot 1)
\end{equation}
for all $x\in\cGh$, $t\in\Mor(\eps,x\tot\ox)$. This antilinear map is well defined: taking $t \in \Mor(\eps,x \tot \ox)$ and $\ttil \in \Mor(\eps,\ox
\tot x)$, normalized in such a way that $(t^* \ot 1)(1 \ot \ttil) = 1$, we define
$$\bigl((\om_{\xi,\eta} \ot \io)(X^x)\bigr)^* := (\om_{(\xi^* \ot
1)t, (1 \ot \eta^*)\vphi(\ttil)} \ot \io)(X^{\ox})$$
for all $\xi\in H_x$ and $\eta\in H_{\vphi(x)}$.

For $\xi\in H_x,\eta\in H_{\vphi(x)}$, we compute
\begin{align*}
\bigl(\bigl((\om_{\xi,\eta} \ot \io)(X^x)\bigr)^*\bigr)^*
&= \bigl((\om_{(\xi^* \ot 1)t, (1 \ot \eta^*)\vphi(\ttil)} \ot
\io)(X^{\ox})\bigr)^* \\ &= (\om_{(((\xi^*\ot1)t)^*\ot
1)\ttil,(1\ot((1\ot\eta^*)\vphi(\ttil))^*\vphi(t))}\ot\io)(X^x)\\
&=(\om_{\xi,\eta}\ot\io)(X^x)\;,
\end{align*} in the last step
using our particular choice of $t$ and $\ttil$.

We also get \[(t^*\ot 1)(X^{\ox}_{23})^*=(\vphi(t)^*\ot 1)X^x_{13}\;.\] Because \[(X^x(X^x)^*)_{13}(\vphi(t)\ot 1)=X^x_{13}X^{\ox}_{23}(t\ot
1)=\vphi(t)\ot 1\] and because, by \eqref{fifu3}, \[(t^*\ot 1)((X^{\ox})^*X^{\ox})_{23}=t^*\ot 1 \; , \] the elements $X^x$ are unitaries.

Since for all $x,y,z\in\cGh$ and $S\in\Mor(z,x\tot y)$, \[(X^x_{13}X^y_{23}(S\ot 1))^*=((\vphi(S)\ot 1)X^z)^*=(X^z)^*(\vphi(S)^*\ot 1)=(S^*\ot
1)(X^y_{23})^*(X^x_{13})^*\] by \eqref{fifu3} and the fact that the $X^x$ are unitary, our involution is anti-multiplicative. We conclude that $\cB$
a *-algebra.

Denote by $\om$ the linear functional $\om : \cB \recht \C$ given by
$\om(1) = 1$ and $(\io \ot \om)(X^x) = 0$ for all $x \neq \eps$. We
show that $\om$ is a faithful state on $\cB$.

Let $x,y \in \cGh$. Take $t \in \Mor(\eps,x \tot \ox)$ such that $t^*t = \dimq(x)$. Take $\ttil \in \Mor(\eps,\ox \tot x)$ such that $(t^* \ot 1)(1
\ot \ttil) = 1$.  Then,
$$(\om_{\mu,\rho} \ot \io)(X^y) (\om_{\xi,\eta} \ot \io)(X^x)^*
= (\om_{\mu \ot (\xi^* \ot 1)t,\rho \ot (1 \ot \eta^*)\vphi(\ttil)}
\ot \io)(X^y_{13} X^{\ox}_{23}) \; .$$ We conclude that
\begin{align*}
\om \bigl( (\om_{\mu,\rho} \ot \io)(X^y)  (\om_{\xi,\eta} \ot
\io)(X^x)^* \bigr) &=\delta_{x,y}
\frac{1}{\dimq(x)} \langle \mu \ot (\xi^* \ot 1) t, t \rangle \;
\langle \vphi(t), \rho \ot (1 \ot \eta^*)\vphi(\ttil) \rangle \\
&=\sde_{x,y}\frac{1}{\dimq(x)}\langle(\xi^*\ot 1)t,(\mu^*\ot 1)t\rangle\langle\rho,\eta\rangle\\
&= \delta_{x,y} \frac{1}{\dimq(x)} \langle Q_x \mu,\xi \rangle \;
\langle \rho,\eta \rangle \; .
\end{align*}
Choose orthonormal bases $(f^x_i)$ for every space $H_{\vphi(x)}$.
Any element $a \in \cB$ admits a unique decomposition
$$a = \sum_{x,i} (\om_{\xi^x_i,f^x_i} \ot \io)(X^x)$$
in terms of vectors $\xi^x_i \in H_x$. Then,
$$\om(aa^*) = \sum_{x,i} \frac{1}{\dimq(x)} \langle Q_x
\xi^x_i,\xi^x_i \rangle \; .$$ It follows that $\om(aa^*) \geq 0$
for all $a$ and that $\om(aa^*) = 0$ if and only if $\xi^x_i = 0$
for all $x$ and $i$, i.e.\ if and only if $a=0$.

The definition of the coaction $\sde : \cB \recht \cB \ot \cA$ is
obvious and it is clear that $\om$ is invariant under $\sde$. It
follows that we can extend $\sde$ to coactions $\sde_r$, resp.\
$\sde_u$ of $(A,\de)$ on $B_r$, resp.\ $B_u$. Moreover, $\om(x) 1 =
(\io \ot h)\sde_r(x)$ for all $x \in B_r$ and analogously for $x \in
B_u$. It follows that $\sde_r$ and $\sde_u$ are ergodic coactions.
Given the unitary elements $X^x$ and Theorem \ref{thm.fullquantum}, it
follows that $\sde_r$ and $\sde_u$ are of full quantum multiplicity.

By definition, the coaction $\sde_r$ on $B_r$ is reduced and the
coaction $\sde_u$ on $B_u$ is universal. Indeed, the canonical
$^*$-subalgebra of $B_u$ generated by the spectral subspaces for
$\sde_u$, is exactly $\cB$.

It remains to show that any reduced, resp.\ universal ergodic coaction of full quantum multiplicity arises as above from a unitary fiber functor. Let
$\sde: B\to B\ot A$ be an ergodic coaction of full quantum multiplicity. Construct for $x \in \cGh$, the unitary elements $X^x \in \B(H_x,\oK_x) \ot
B$ as in Section \ref{sec.spectral}. Define $H_{\vphi(x)} := \oK_x$. Let $S \in \Mor(x_1 \tot \cdots \tot x_r,y_1 \tot \cdots \tot y_k)$. The element
$$X^{y_1}_{1,k+1} \cdots X^{y_k}_{k,k+1} (S \ot 1)(X^{x_1}_{1,r+1}
\cdots X^{x_r}_{r,r+1})^*$$
is invariant under $\io \ot \sde$. So, we can define $\vphi(S)$ by the formula
$$X^{y_1}_{1,k+1} \cdots X^{y_k}_{k,k+1} (S \ot 1) = (\vphi(S) \ot 1) X^{x_1}_{1,r+1}
\cdots X^{x_r}_{r,r+1}$$ for all $S \in \Mor(x_1 \tot \cdots \tot x_r,y_1 \tot \cdots \tot y_k)$. It is clear that $\vphi$ is a unitary fiber functor
on $\cG$.

Denote by $\cB$ the $^*$-subalgebra of $B$ generated by the spectral
subspaces of $\sde$ as defined in Section \ref{sec.spectral}. By
definition, $\cB$ is generated by the coefficients of the $X^x$. In
order to show that $\cB$ is isomorphic to the $^*$-algebra defined
by the unitary fiber functor $\vphi$, it suffices to show that the
coefficients of the $X^x$ form a linear basis of $\cB$. But this
follows immediately from \eqref{eq.invariant-state}.
\end{proof}

\begin{definition}
Two unitary fiber functors $\vphi$ and $\psi$ on a compact quantum group $\cG$ are said to be \emph{isomorphic} if there exist unitaries $u_x \in
\B(H_{\vphi(x)},H_{\psi(x)})$ satisfying
$$\psi(S) = (u_{y_1} \ot \cdots \ot u_{y_k}) \vphi(S) (u_{x_1}^* \ot \cdots \ot u_{x_r}^*)$$
for all $S \in \Mor(x_1 \tot \cdots \tot x_r,y_1 \tot \cdots \tot y_k)$.
\end{definition}

\begin{proposition}
Let $\vphi$ and $\psi$ be unitary fiber functors on $\cG$ and denote
by $\sde_\vphi$, resp.\ $\sde_\psi$ the associated coactions on
$\cB_\vphi$, resp.\ $\cB_\psi$. Then the following statements are
equivalent.
\begin{itemize}
\item The fiber functors $\vphi$ and $\psi$ are isomorphic.
\item There exists a $^*$-isomorphism $\pi : \cB_\vphi \recht
  \cB_\psi$ satisfying $(\pi \ot \io)\sde_{\vphi} = \sde_\psi \pi$.
\end{itemize}
\end{proposition}

\begin{proof}
Straightforward.
\end{proof}

The following proposition follows immediately from the Tannaka-Krein
reconstruction theorem, but we give a detailed statement for
clarity. Its proof is completely analogous to the proof of Theorem \ref{thm.main}.

\begin{proposition} \label{prop.tannaka-krein}
Let $\vphi$ be a unitary fiber functor on a compact quantum group $\cG = (A,\de)$.
\begin{itemize}
\item There exists a unique universal compact quantum group $(A_2,\de_2)$ with underlying Hopf $^*$-algebra \linebreak $(\cA_2,\de_2)$ and unitary
  corepresentations $U^{\vphi(x)} \in \B(H_{\vphi(x)}) \ot \cA_2$
  satisfying
\begin{enumerate}
\item $U^{\vphi(y)}_{13} U^{\vphi(z)}_{23} (\vphi(S) \ot 1) =
  (\vphi(S) \ot 1) U^{\vphi(x)} \quad\text{for all}\quad S \in
  \Mor(x,y \tot z)$,
\item the matrix coefficients of the $U^{\vphi(x)}$ form a linear
  basis of $\cA_2$.
\end{enumerate}
\item $\{U^{\vphi(x)} \mid x \in \cG\}$ is a complete set of
  irreducible corepresentations of $(A_2,\de_2)$ and $\vphi$ is a
  monoidal equivalence of compact quantum groups.
\end{itemize}
\end{proposition}

The following result is then a corollary of Theorem \ref{thm.main}.

\begin{proposition} \label{prop.pair}
Consider two compact quantum groups $\cG=(A,\de)$ and
$\cG_2=(A_2,\de_2)$. Let $\vphi:\cG\to\cG_2$ be a monoidal
equivalence. In particular, $\vphi$ is a unitary fiber functor on $\cG$.

Denote by $B_r$, resp.\ $B_u$ the \cst-algebras associated to $\vphi$
as in Theorem \ref{thm.main}, with dense $^*$-subalgebra $\cB$.
Denote by $\sde$ the
corresponding coaction of $(\cA,\de)$ on $\cB$. Denote by $X^x \in
\B(H_x,H_{\vphi(x)}) \ot \cB$ the unitaries generating $\cB$.
\begin{itemize}
\item There is a unique coaction $\sde_2 : \cB \recht \cA_2 \ot \cB$
  satisfying $(\io \ot \sde_2)(X^x) = U^{\vphi(x)}_{13} X^x_{23}$ for
  all $x \in \cGh$. The coaction $\sde_2$ commutes with $\sde$ and extends to $B_r$, resp.\
  $B_u$, yielding a reduced, resp.\ universal, ergodic coaction of
  full quantum multiplicity.
\item Every pair of commuting reduced, resp.\ universal, ergodic
  coactions of full quantum multiplicity arises in this way from a
  monoidal equivalence.
\end{itemize}
\end{proposition}

\begin{proof}
Given the monoidal equivalence $\vphi$, it is obvious to construct
the coaction $\sde_2$.

It remains to show the second statement. Let $\sde : B \recht B \ot A$
and $\sde_2 : B \recht A_2 \ot B$ be commuting ergodic coactions of
full quantum multiplicity. Denote by $\cB$ the unital $^*$-subalgebra
of $B$ generated by the spectral subspaces of $\sde$.
Using Theorem \ref{thm.main},
we get a unitary fiber functor $\vphi$ on $\cG$ and we may assume that
$\cB$ and $\sde$ are constructed from $\vphi$ as in Theorem
\ref{thm.main}. In particular, $\cB$ is generated by the coefficients
of $X^x \in \B(H_x,H_{\vphi(x)}) \ot \cB$.

Because $\sde$ and $\sde_2$ commute, the element
$(\io\ot\sde_2)(X^x)(X^x)^*_{13}$ is invariant under
$(\io\ot\io\ot\sde)$. Since $\sde$ is ergodic and $X^x$ unitary,
we get a unitary element $U^{\vphi(x)}\in
B(H_{\vphi(x)})\ot A_2$ such that $(\io \ot \sde_2)(X^x) =
U^{\vphi(x)}_{12} X^x_{13}$. Because $\sde_2$ is a coaction, we easily
compute that $U^{\vphi(x)}$ is a unitary corepresentation of $\cG_2$.

It remains to show that $\{U^{\vphi(x)} \mid x \in \cGh \}$ is a complete set of irreducible unitary corepresentations of $\cG_2$ and that $\vphi$ is
a monoidal equivalence.

Assume that $S\in\Mor(\vphi(x),\vphi(y))$. The element $(X^y)^*(S\ot 1)X^x\in \B(H_x,H_y)\ot B$ is
invariant under $\io\ot\sde_2$, so it has the form $T\ot 1$, with
$T\in \B(H_x,H_y)$. It follows that $T \in \Mor(x,y) = \sde_{x,y} \C$
and hence, $S \in \sde_{x,y} \C$. So, the $U^{\vphi(x)}$ are mutually
inequivalent irreducible corepresentations of $\cG_2$.

In order to show that the set $\{U^{\vphi(x)} \mid x \in \cGh\}$
exhausts all irreducible corepresentations of $\cG_2$, it suffices to
show, for all $a \in \cA_2$, that $(\io \ot h_2)((1 \ot
a)U^{\vphi(x)}) = 0$ for all $x \in \cGh$, implies $a =0$. But, given
the formula for $\sde_2$, we get $(h_2 \ot \io)((a \ot 1)\sde_2(x)) =
0$ for all $x \in B$. Since $\sde_2$ is of full quantum multiplicity,
this implies that $a=0$.

It remains to show that $\vphi$ is a monoidal equivalence. For this, it suffices to show that $\Mor(\vphi(x),\vphi(y) \tot \vphi(z)) = \vphi(\Mor(x,y
\tot z))$. If $S \in \Mor(x,y \tot z)$, we use the multiplicativity of $\sde_2$ to obtain
\begin{align*}
(\vphi(S) \ot 1 \ot 1) U^{\vphi(x)}_{12} X^x_{13}
&= (\io \ot \io \ot \sde_2)( (\vphi(S) \ot 1) X^x)
= (\io \ot \io \ot \sde_2)( X^y_{13} X^z_{23} (S \ot 1) )
\\ &= U^{\vphi(y)}_{13} X^y_{14} \; U^{\vphi(z)}_{23} X^z_{24} \; (S \ot 1
\ot 1)
= U^{\vphi(y)}_{13} U^{\vphi(z)}_{23} (\vphi(S) \ot 1 \ot 1) X^x_{13} \; .
\end{align*}
It follows that $S \in \Mor(\vphi(x),\vphi(y) \tot \vphi(z))$. The converse inclusion is shown analogously.
\end{proof}

\section{Unitary fiber functors preserving the dimension} \label{sec.cocycle}

We study in this section unitary fiber functors $\vphi$ on a compact quantum group $\cG$ preserving the dimension, i.e.\ satisfying $\dim
H_{\vphi(x)} = \dim H_x$ for all $x \in \cGh$. Taking into account Theorem \ref{thm.main}, this comes down to the study of ergodic coactions of full
quantum multiplicity satisfying $\mult(x) = \dim(x)$ for all $x \in \cGh$.

We establish a relation between unitary fiber functors preserving the dimension (up to isomorphism) and the $2$-cohomology of the dual, discrete
quantum group $(\Ah,\deh)$. The following definition is due to Landstad \cite{landstad2} and Wassermann \cite{wasser2}, who consider it for the dual
of a compact group.

\begin{definition}
A unitary element $\Om \in \M(\Ah \ot \Ah)$ is said to be a \emph{$2$-cocycle} if it satisfies
\begin{equation}\label{eq.cocycle}
(\deh \ot \io)(\Om) (\Om \ot 1) =
(\io \ot \deh)(\Om)(1 \ot \Om)\;.
\end{equation}
Two $2$-cocycles $\Om_1$ and $\Om_2$ are said to differ by a coboundary if there exists a unitary $u \in \M(\Ah)$ such that $\Om_2 = \deh(u) \Om (u^*
\ot u^*)$. We denote this relation by $\Om_1 \sim \Om_2$ and observe that $\sim$ is an equivalence relation on the set of $2$-cocycles.
\end{definition}

\begin{remark}
In the quantum setting, there is no reason that the product of two $2$-cocycles is again a $2$-cocycle. So, although we could define the
$2$-cohomology of $(\Ah,\deh)$ as the set of equivalence classes of $2$-cocycles, this set has no natural group structure.
\end{remark}

\begin{notation}
We denote by $p_x, x \in \cGh$, the minimal central projections of $\Ah = \oplus_x \B(H_x)$.
\end{notation}

\begin{remark}
Up to coboundary, we can and will assume that a unitary $2$-cocycle is \emph{normalized}, i.e.
$$(p_\eps \ot 1) \Om = p_\eps \ot 1 \quad\text{and}\quad (1 \ot p_\eps)\Om = 1 \ot p_\eps \; .$$
\end{remark}

Let $\Om$ be a normalized unitary $2$-cocycle on $(\Ah,\deh)$, the dual of $\cG = (A,\de)$. Denote
$$\Om_{(2)} := (\deh \ot \io)(\Om) (\Om \ot 1) = (\io \ot \deh)(\Om) (1 \ot \Om) \; .$$
It follows from Remark \ref{rem.reduction} that there is a unique unitary fiber functor $\vphi_\Om$ on $\cG$ satisfying
$$H_{\vphi_\Om(x)} = H_x \; , \quad \vphi_\Om(S) = \Om^* S \; , \quad \vphi_\Om(T) = \Om_{(2)}^* T \; ,$$
for all $S \in \Mor(x,y \tot z)$ and $T \in \Mor(a,x \tot y \tot z)$. Observe that we implicitly used that $\B(H_y \ot H_z)$ is an ideal in $\M(\Ah
\ot \Ah)$ and hence, $\Om^* S$ is a well defined element of $\B(H_x,H_y \ot H_z)$ whenever $S \in \Mor(x,y \tot z)$.

From Proposition \ref{prop.tannaka-krein}, we get a compact quantum group $(A_\Om,\de_\Om)$, whose dual $(\Ah_\Om,\deh_\Om)$ is given by
$$\Ah_\Om = \oplus_x \B(H_x) = \Ah \quad\text{and}\quad \deh_\Om(a) \vphi_\Om(S) = \vphi_\Om(S) a \quad\text{for all}\quad a \in \B(H_x), S \in
\Mor(x, y \tot z) \; .$$ Also, $\vphi_\Om$ becomes a monoidal equivalence between $(A,\de)$ and $(A_\Om,\de_\Om)$. Observe that
$$\deh_\Om(a) = \Om^* \deh(a) \Om \quad\text{for all}\quad a \in \Ah_\Om = \Ah \; .$$

\begin{proposition}
Let $\vphi$ be a unitary fiber functor on a compact quantum group $\cG = (A,\de)$ such that $\dim H_{\vphi(x)} = \dim H_x$ for all $x \in \cGh$. Then
there exists a normalized unitary $2$-cocycle $\Om$ on $(\Ah,\deh)$, uniquely determined up to coboundary, such that $\vphi$ is isomorphic with
$\vphi_\Om$.
\end{proposition}
\begin{proof}
Denote by $\sde : B \recht B \ot A$ the reduced ergodic coaction associated with $\vphi$ by Theorem \ref{thm.main}. Consider the generating unitaries
$X^x \in \B(H_x,H_{\vphi(x)}) \ot B$ satisfying $(\io \ot \sde)(X^x) = X^x_{12} U^x_{13}$ for all $x \in \cGh$.

Since $\dim H_{\vphi(x)} = \dim H_x$, we can take unitary elements $u_x : H_{\vphi(x)} \recht H_x$. Take $u_\eps = 1$. Define $Y^x = (u_x \ot 1)X^x$
and consider $Y:= \oplus_x Y^x \in \M(\Ah \ot B)$. Because the element $(\deh \ot \io)(Y) Y^*_{23} Y^*_{13}$ is invariant under $(\io \ot \io \ot
\sde)$, we find a unitary element $\Om \in \M(\Ah \ot \Ah)$ such that
$$(\deh \ot \io)(Y) = (\Om \ot 1)Y_{13}Y_{23} \; .$$
Applying $\deh \ot \io \ot \io$ and $\io \ot \deh \ot \io$ to this equality, we obtain that $\Om$ is a unitary $2$-cocycle on $(\Ah,\deh)$.

It remains to show that $\vphi$ and $\vphi_\Om$ are isomorphic. Let $S \in \Mor(x,y \tot z)$. Then,
\begin{align*}
(S \ot 1)Y^x &= (\deh \ot \io)(Y)(S \ot 1) = (\Om \ot 1)Y_{13}Y_{23} (S \ot 1) = (\Om(u_y \ot u_z) \ot 1) X^y_{13} X^z_{23} (S \ot 1) \\ &= (\Om(u_y
\ot u_z) \ot 1) (\vphi(S) \ot 1) X^x = (\Om(u_y \ot u_z) \vphi(S) u_x^* \ot 1) Y^x \; .
\end{align*}
Hence, $\vphi_\Om(S) = \Om^* S = (u_y \ot u_z) \vphi(S) u_x^*$ for all $S \in \Mor(x,y\tot z)$.

It is obvious that $\vphi_{\Om_1}$ is isomorphic with $\vphi_{\Om_2}$ if and only if the $2$-cocycles $\Om_1$ and $\Om_2$ differ by a coboundary.
\end{proof}

Fix a normalized unitary $2$-cocycle $\Om$ on $(\Ah,\deh)$ and consider the unitary fiber functor $\vphi_\Om$. Theorem \ref{thm.main} yields
\cst-algebras $B^\Om_r$ and $B^\Om_u$ with ergodic coactions $\sde_r$ and $\sde_u$ of full quantum multiplicity. It is, of course, possible to
describe these \cst-algebras directly in terms of $\Om$: they correspond to the $\Om$-twisted group \cst-algebras of $(\Ah,\deh)$. Before we can
prove such a statement, we have to introduce a few notations and a bit of terminology. Such $\Om$-twisted group \cst-algebras have been studied by
Landstad \cite{landstad2} and Wassermann \cite{wasser2} when $\Om$ is a unitary $2$-cocycle on the dual of a compact group.

\begin{notation} \label{not.multun}
Denote by $H$ the $L^2$-space of the Haar state $h$ on $(A,\de)$. We consider $A$ and $\Ah$ as being represented on $H$. Denote by $V$ the right
regular corepresentation of $(A,\de)$. The unitary $V$ is
multiplicative in the sense of \cite{BS} and belongs to $\M(\Ah \ot A)$. There exists a canonical unitary $u \in \B(H)$
with $u^2 = 1$ such that $\Vh := (u \ot 1) \Sigma V \Sigma (u \ot 1)$ is a multiplicative unitary. We know that $u A u \subset A'$, $u \Ah u
\subset \Ah'$ and
$$\deh(x) = V^* (1 \ot x)V = \Vh(x \ot 1) \Vh^*
\quad\text{and}\quad (\io \ot \de)(V) = V_{12} V_{13} \; , (\io \ot \deh)(\Vh) = \Vh_{12} \Vh_{13} \; .$$ The unitaries $V$ and $\Vh$ satisfy
the pentagonal equation
\[V_{12}V_{13}V_{23}=V_{23}V_{12}\quad\text{and}\quad \Vh_{12}\Vh_{13}\Vh_{23}=\Vh_{23}\Vh_{12}\]
\end{notation}

\begin{definition}
An $\Om$-corepresentation of $(\Ah,\deh)$ on a Hilbert space $K$ is a unitary $X \in \M(\Ah \ot \cK(K))$ satisfying
$$(\deh \ot \io)(X) = (\Om \ot 1) X_{13} X_{23} \; .$$
\end{definition}

The following lemma can be checked immediately using the formulas in Notation \ref{not.multun} and in particular, the commutation relation $V_{12}
\Vh_{23} = \Vh_{23} V_{12}$.
\begin{lemma}
Denoting $\Omtil := (1 \ot u) \Sigma \Om^* \Sigma (1 \ot u)$, the
unitary $\Omtil V \in \M(\Ah \ot \cK(H))$ is an
$\Om$-corepresentation. It is called the
\emph{right regular $\Om$-corepresentation}.
\end{lemma}

The next lemma is crucial to define the twisted quantum group \cst-algebras.

\begin{lemma}
Let $X$ be an $\Om$-corepresentation of $(\Ah,\deh)$ on $K$. Then, $[(\mu \ot \io)(X) \mid \mu \in \Ah^*]$ is a unital \cst-algebra.
\end{lemma}
\begin{proof}
Write $B:= [(\mu \ot \io)(X) \mid \mu \in \Ah^*]$. From the defining relation for an $\Om$-corepresentation, it follows that $B$ is an algebra acting
non-degenerately on $K$. Since $(\deh \ot \io)(X) X^*_{23} = (\Om \ot 1)X_{13}$, we have
\begin{align*}
B &= [(\mu \ot \eta \ot \io)\bigl( (\deh \ot \io)(X) X^*_{23} \bigr) \mid \mu,\eta \in \Ah^*] = [(\mu \ot \eta \ot \io)\bigl( \Vh_{12} X_{13}
\Vh^*_{12} X^*_{23} \bigr) \mid \mu,\eta \in \Ah^*] \\&=[(\mu \ot \eta \ot \io)\bigl( X_{13} \Vh^*_{12} X^*_{23} \bigr) \mid \mu,\eta \in \Ah^*] =
[(\mu \ot \eta \ot \io)\bigl(X_{13} \bigl((1 \ot \cK(H))\Vh^*(\cK(H) \ot 1)\bigr)_{12} X^*_{23} \bigr) \mid \mu,\eta \in \Ah^*] \\ &= [(\mu \ot \eta
\ot \io)\bigl(X_{13} (\cK(H) \ot \cK(H) \ot 1) X^*_{23}
\bigr) \mid \mu,\eta \in \Ah^*] \\
&=[(\mu \ot \io)(X) (\eta \ot \io)(X)^* \mid \mu,\eta \in \Ah^*] = [B B^*] \; .
\end{align*}
Here we used the regularity of the multiplicative unitary $\Vh$. It follows that $B$ is a \cst-algebra. Since $\Om$ is normalized, we have $(p_\eps
\ot 1) X = p_\eps \ot 1$ and hence, $B$ is unital.
\end{proof}

\begin{definition}
We define the unital \cst-algebras
$$\cstr(\cG,\Om):= [(\mu \ot \io)(\Omtil V) \mid \mu \in \Ah^*] \quad\text{and}\quad
\cstu(\cG,\Om):=[(\mu \ot \io)(X) \mid \mu \in \Ah^*] \; ,$$ where $X$ denotes a universal $\Om$-corepresentation.
\end{definition}

Remark that an $\Om$-corepresentation $X$ on $K$ is said to be \emph{universal} if for any $\Om$-corepresentation $Y$ on $K_1$, there exists an isometry $v
: K_1 \recht K$ such that $X(1 \ot v) = (1 \ot v)Y$. It is clear that a universal $\Om$-corepresentation exists and that the \cst-algebra
$\cstu(\cG,\Om)$ is well defined up to an isomorphism.

\begin{proposition}
Denote by $\cB$ the unital $^*$-algebra associated by Theorem \ref{thm.main} with the unitary fiber functor $\vphi_\Om$. Consider the unitaries $X^x
\in \B(H_x) \ot \cB$ generating $\cB$. Denote by $B^\Om_r$ and $B^\Om_u$ the associated reduced and universal \cst-algebra.

For any $\Om$-corepresentation $X$ of $(\Ah,\deh)$ on a Hilbert space $K$, we obtain a representation $\pi$ of $\cB$ on $K$ given by $(\io \ot
\pi)(X^x) = (p_x \ot 1)X$ for all $x \in \cGh$. Taking $X = \Omtil V$, we get an isomorphism $B^\Om_r \cong \cstr(\cG,\Om)$. Taking $X$ to be a
universal $\Om$-corepresentation, we get an isomorphism $B^\Om_u \cong \cstu(\cG,\Om)$.
\end{proposition}

\begin{proof}
It is immediate that the formula $(\io \ot \pi)(X^x) = (p_x \ot 1)X$ defines a one-to-one correspondence between representations of $\cB$ and
$\Om$-corepresentations of $(\Ah,\deh)$. This already shows the isomorphism $B^\Om_u \cong \cstu(\cG,\Om)$.

Consider $\pi : \cB \recht \B(H)$ given by $(\io \ot \pi)(X^x) = (p_x \ot 1)\Omtil V$. Denote by $\om$ the unique invariant state on $\cB$. To prove
the isomorphism $B^\Om_r \cong \cstr(\cG,\Om)$, it suffices to define a faithful state $\om_1$ on $\cstr(\cG,\Om)$ such that $\om_1 \pi = \om$.
Define $\al : \cstr(\cG,\Om) \recht \M(\cK(H) \ot A) : \al(a) = V(a \ot 1)V^*$. One verifies that $(\io \ot \al)(\Omtil V) = (\Omtil V)_{12} V_{13}$.
Hence, $\al : \cstr(\cG,\Om) \recht \cstr(\cG,\Om) \ot A$ is a coaction satisfying $\al \pi = (\pi \ot \io)\sde$. It follows that $(\io \ot h)\al(a)
\in \C 1$ for all $a \in \cstr(\cG,\Om)$. So, we can define $\om_1$ by the formula $\om_1(a)1 = (\io \ot h)\al(a)$. Clearly, $\om_1 \pi = \om$.
\end{proof}

\section{Monoidal equivalence for $A_o(F)$} \label{sec.ao}

Recall the following definition of the compact quantum group $A_o(F)$ \cite{VDW}.

\begin{definition}
For all $n\in\mathbb{N}$ and $F\in \GL(n,C)$ with $F\overline{F} = c 1 \in \R 1$, $A_0(F)$ is defined as the universal compact quantum group
generated by the coefficients of the corepresentation $U \in M_n(\C)
\ot A_o(F)$ with
relations
$$U \;\;\text{is unitary and}\;\; U = (F \ot 1)\overline{U} (F^{-1}
\ot 1) \; .$$
Then,
$(A_o(F),U)$ is a compact (matrix) quantum group.
\end{definition}

When the matrix $F$ has dimension $2$, we precisely obtain the quantized versions of the classical Lie group $\SU(2)$, as considered by Woronowicz.

\begin{definition}(\cite{wor})
Let $q \in [-1,1] \setminus \{0\}$. We define $\SU_q(2)$ to be the universal \cst-algebra generated by $2$ elements $\alpha, \gamma$ such that
$U=\Bigl(\begin{smallmatrix} \alpha & -q\gamma^* \\ \gamma &\alpha^*
\end{smallmatrix}\Bigr)$ is a unitary corepresentation. Then $(\SU_q(2),U)$ is a compact (matrix)
quantum group.
\end{definition}

Observe that $\SU_q(2) \cong A_o\Bigl(\begin{smallmatrix} 0 & 1 \\ -q^{-1}
  & 0 \end{smallmatrix} \Bigr)$ and that $\dimq(U) = \bigl| q + \frac{1}{q} \bigr|$.

Banica \cite{banica1} has shown that the irreducible corepresentations of $A_o(F)$ can be labelled by $\N$, in such a way that the fusion rules are
identical to the fusion rules of the compact Lie group $\SU(2)$. In particular, $\Mor(\eps,U \tot U)$ is one-dimensional and generated by
$$t_F := \sum_{i} e_i \ot F e_i \; ,$$
where $(e_i)$ is the standard basis of $\C^n$.

\begin{theorem} \label{thm.fiberao}
Let $F\in \GL(n,\C)$ with $F\overline{F}=c 1$ and $c \in \R$. Consider $\cG = A_o(F)$.
\begin{itemize}
\item Take $F_1 \in \GL(n_1,\C)$ satisfying $F_1 \overline{F}_1 =
  c_1 1$ and $\frac{c}{\Tr(F^*F)}=\frac{c_1}{\Tr(F_1^*F_1)}$.
 There exists a unitary fiber functor $\vphi_{F_1}$ on $\cG$,
  uniquely determined up to isomorphism, such that
\begin{equation}\label{eq.fiberao}
\vphi\bigl( \frac{1}{\sqrt{\Tr(F^*F)}} t_F \bigr) = \frac{1}{\sqrt{\Tr(F_1^*F_1)}} t_{F_1} \; .
\end{equation}
\item Every unitary fiber functor $\vphi$ on $\cG$ is isomorphic with
  one of the form $\vphi_{F_1}$. Moreover, $\vphi_{F_1}$ is isomorphic
  with $\vphi_{F_2}$ if and only if $n_1=n_2$ and there exists a
  unitary $v \in \U(n_1)$ and a $\lambda \neq 0$ such that $F_2 =
  \lambda v F_1 v^t$.
\end{itemize}
\end{theorem}

\begin{proof}
Take a parameter $\be \in \R \setminus
\{0\}$. Consider the universal graded \cst-algebra $(A(n,m))_{n,m
\in \N}$ satisfying $A_{n,m} = \{0\}$ if $n -m$ is odd and generated
by elements $v(r,s) \in A(r+s,r+s+2)$ with relations (denoting $1_n$
the unit of the \cst-algebra $A(n) := A(n,n)$)
\begin{align*}
v(r,s)^* v(r,s) &= 1_{r+s} \\
v(r,s+1)^* v(r+1,s) &= \be \; 1_{r+s+1} \\
v(r,k+l+2) v(r+k,l) &= v(r+k+2,l) v(r,k+l) \\
v(r,k+l+2)^* v(r+k+2,l) &= v(r+k,l) v(r,k+l)^*
\end{align*}

Take $F \in \GL(n,\C)$ with $F \overline{F} = c 1$ and $c \in \R$. Put $\be =
\frac{c}{\Tr(F^* F)}$. Let $\cG=(A_o(F),\de)$ and denote by $U^n$
the $n$-fold tensor product of the fundamental corepresentation, with
the convention that $U^0=\eps$. Take the isometric $t \in \Mor(\eps,U^2)$ given by
$t = \frac{1}{\sqrt{\Tr(F^*F)}} t_F$. Then, $(t^* \ot 1)(1 \ot t) = \be \; 1$.

We get a natural $^*$-homomorphism
$$\pi : (A(n,m))_{n,m \in \N} \recht (\Mor(U^n,U^m))_{n,m \in \N}$$
given by $\pi(v(r,s)) = 1_r \ot t \ot 1_s$. Because of Proposition 1
in \cite{banica2}, $\pi$ is surjective. It follows
from the comments after Th{\'e}or{\`e}me 4 in \cite{banica1}, that $\pi$ is
faithful on $A(n,n)$ for all $n$. But then, $\pi$ is faithful on
$A(n,m)$ because
\[\pi(T)=0\Leftrightarrow\pi(T^*T)=0\Leftrightarrow T^*T=0\Leftrightarrow
T=0\;.\]
We conclude that $\pi$ is a $^*$-isomorphism.

Take $F_1 \in \GL(n_1,\C)$ satisfying $F_1 \overline{F}_1 =
  c_1 1$ and $\frac{c}{\Tr(F^*F)}=\frac{c_1}{\Tr(F_1^*F_1)}$. Write
  $K=\C^{n_1}$ and denote by $K^n$ the $n$-fold tensor product of $K$,
  with the convention that $K^0 = \C$.
From the
  preceding discussion, we obtain a faithful $^*$-homomorphism
$$\pi : (\Mor(U^n,U^m))_{n,m \in \N} \recht (\B(K^n,K^m))_{n,m \in
  \N}$$ satisfying $\pi(t) = \frac{1}{\sqrt{\Tr(F_1^*F_1)}}
  t_{F_1}$, $\pi(1) = 1$ and $\pi(1 \ot T \ot 1) = 1 \ot \pi(T) \ot 1$ for all $T$.

We choose a concrete identification $\cGh = \N$ as follows. We define $P_n \in \Mor(U^n,U^n)$ as the unique projection satisfying $P_n T = 0$ for all
$r < n$ and all $T \in \Mor(U^r,U^n)$. We define $U_n$ as the restriction of $U^n$ to the image of $P_n$. We then identify
$$\Mor(n_1 \tot \cdots \tot n_r, m_1 \tot \cdots \tot m_k) = (P_{m_1} \ot
\cdots \ot P_{m_k}) \Mor(U^{n_1 + \cdots + n_r},U^{m_1 + \cdots +
  m_k}) (P_{n_1} \ot \cdots \ot P_{n_r}) \; .$$
Define $H_{\vphi(n)} := \pi(P_n)K^n$ and define $\vphi(S)$ by restricting $\pi$ to $\Mor(n_1 \tot \cdots \tot n_r, m_1 \tot \cdots \tot m_k)$. It is
now obvious that $\vphi$ is a unitary fiber functor on $\cG$.

Suppose conversely that $\vphi$ is a unitary fiber functor on $A_o(F)$. We continue to use the concrete identification $\cGh = \N$ introduced above.
Up to isomorphism, we may assume that $H_{\vphi(1)} = \C^{n_1}$ and we denote $K = \C^{n_1}$. We define the $^*$-homomorphism $\pi :
(\Mor(U^n,U^m))_{n,m} \recht (\B(K^n,K^m))_{n,m}$ by restricting $\vphi$. Define the matrix $F_1$ such that $\vphi(t_F) = t_{F_1}$. Then, $F_1
\overline{F}_1 = c_1 1$ with $c_1 = c$ and $\Tr(F_1^*F_1) = \Tr(F^* F)$. Since $t$ generates the graded \cst-algebra $(\Mor(U^n,U^m))_{n,m}$, $\pi$
coincides with the $^*$-homomorphism constructed in the first part of the proof starting with $F_1$. Denoting by $T_n \in \Mor(n,1 \tot \cdots \tot
1)$ the embedding, we get unitary operators $\vphi(T_n) : H_{\vphi(n)} \recht \pi(P_n)K^n$ who implement the isomorphism between $\vphi$ and
$\vphi_{F_1}$.
\end{proof}

\begin{corollary} \label{cor.moneqao}
Let $F\in \GL(n,\C)$ with $F\overline{F}=c 1$ and $c \in \R$.
Consider $\cG = A_o(F)$. A compact quantum group $\cG_1$ is monoidally equivalent with $\cG$ if and
only if there exists $F_1 \in \GL(n_1,\C)$ satisfying $F_1 \overline{F}_1 =
  c_1 1$ and $\frac{c}{\Tr(F^*F)}=\frac{c_1}{\Tr(F_1^*F_1)}$ such that
  $\cG_1 \cong A_o(F_1)$.
\end{corollary}

\begin{proof}
The unitary fiber functor $\vphi_{F_1}$ constructed in Theorem
\ref{thm.fiberao} yields a monoidal equivalence $A_o(F) \meq
A_o(F_1)$. Since these fiber functors $\vphi_{F_1}$ are, up to
isomorphism, the only unitary
fiber functors on $A_o(F)$, we are done.
\end{proof}

So, we exactly know when the compact quantum groups $A_o(F_1)$ and $A_o(F_2)$ are monoidally equivalent. If this is the case, Proposition
\ref{prop.pair} provides us with a universal \cst-algebra $B_u$ and a pair of ergodic coactions of full quantum multiplicity. It is possible to give
an intrinsic description of this \cst-algebra $B_u$.

\begin{theorem} \label{thm.aoff}
Let $F_i \in \GL(n_i,\C)$ by such that $F_i \overline{F}_i = c_i 1$ for $c_i \in \R$ ($i=1,2$). Assume that $c_1 = c_2$ and $\Tr(F_1^* F_1) =
\Tr(F_2^* F_2)$.
\begin{itemize}
\item Denote by $A_o(F_1,F_2)$ the universal unital \cst-algebra
  generated by the coefficients of
$$Y \in M_{n_2,n_1}(\C) \ot A_o(F_1,F_2) \quad\text{with relations}\quad
Y \;\;\text{unitary}\quad\text{and}\quad Y = (F_2 \ot 1)\overline{Y} (F_1^{-1} \ot 1) \; .$$ Then, $A_o(F_1,F_2) \neq 0$ and there exists a unique
pair of commuting universal ergodic coactions of full quantum multiplicity, $\sde_1$ of $A_o(F_1)$ and $\sde_2$ of $A_o(F_2)$, such that
$$(\io \ot \sde_1)(Y) = Y_{12} (U_1)_{13} \quad\text{and}\quad (\io
\ot \sde_2)(Y) = (U_2)_{12} Y_{13} \; .$$ Here, $U_i$ denotes the
fundamental corepresentation of $A_o(F_i)$.
\item $(A_o(F_1,F_2),\sde_1,\sde_2)$ is isomorphic with the \cst-algebra $B_u$ and the coactions thereon given by Proposition \ref{prop.pair} and the
monoidal equivalence $A_o(F_1) \meq A_o(F_2)$ of Corollary \ref{cor.moneqao}.
\item The multiplicity of the fundamental corepresentation $U_1$ in the coaction $\sde_1$ equals $n_2$.
\end{itemize}
\end{theorem}

Remark that the condition on the matrices $F_1$ and $F_2$ is not really less general then the condition in Theorem \ref{thm.fiberao}, but just a
normalization: if $\frac{c_1}{\Tr(F_1^*F_1)}=\frac{c_2}{\Tr(F_2^*F_2)}$, we multiply $F_2$ by a scalar and obtain $c_1 = c_2$ and $\Tr(F_1^* F_1) =
\Tr(F_2^* F_2)$.

\begin{proof}
Take $F_i$ as in the statement of the theorem and denote by $\vphi$ the unitary fiber functor on $\cG_1:=A_o(F_1)$ given by Theorem \ref{thm.fiberao}
and \eqref{eq.fiberao}. We continue to use the identification of $\cGh_1$ with $\N$. Theorem \ref{thm.main} provides us with a $^*$-algebra $\cB$
generated by the coefficients of unitary operators $X^n \in \B(H_n,H_{\vphi(n)}) \ot \cB$, $n \in \N$. Define $X:=X^1$. Since every element of
$\cGh_1$ appears in a tensor power of the fundamental corepresentation and since $\overline{1} = 1$, it follows that the coefficients of $X$ generate
$\cB$ as an algebra. Moreover \eqref{involution} precisely says that $X=(F_2\ot 1)\overline{X}(F_1^{-1}\ot 1)$. It follows that $A_o(F_1,F_2) \neq
0$. Denoting by $\cC$ the unital $^*$-subalgebra of $A_o(F_1,F_2)$ generated by the coefficients of $Y$, we get a surjective $^*$-homomorphism $\rho
: \cC \recht \cB$ satisfying $(\io \ot \rho)(Y) = X$. It remains to show that $\rho$ is a $^*$-isomorphism.

Denote by $U$ the fundamental corepresentation of $A_o(F_1)$ on $H=\C^{n_1}$ and by $U^n$ its $n$-th tensor power, on $H^n$. As in the proof of
Theorem \ref{thm.fiberao}, we denote by $P_n \in \Mor(U^n,U^n)$ the projection onto the irreducible corepresentation $U_n$. Denote $K = \C^{n_2}$ and
denote by $K^n$ the $n$-th tensor power of $K$. Recall that we constructed a faithful $^*$-homomorphism $\pi : (\Mor(U^n,U^m))_{n,m} \recht
(\B(K^n,K^m))_{n,m}$ and that $\vphi$ is defined by restricting $\pi$ to the relevant subspaces. The graded \cst-algebra $(\Mor(U^n,U^m))_{n,m}$ is
generated by the elements $t \in \Mor(\eps,U^2)$ and $1_r \ot t \ot 1_s$.

Put $Y^n=Y_{1,n+1}\cdots Y_{n,n+1} \in \B(H^n,K^n) \ot A_o(F_1,F_2)$. Since $Y^2(t \ot 1) = \pi(t) \ot 1$, it follows that $Y^m(S\ot 1)=(\pi(S)\ot
1)Y^n$ for all $S\in\Mor(U^n,U^m)$. Define a linear map $\gamma : \cB \recht \cC$ by the formula $(\io\ot\gamma)(X^n)=Y^n(P_n \ot 1)=(\pi(P_n) \ot
1)Y^n$. We claim that $\gamma$ is multiplicative. Take $a,n,m\in\N$ and $S\in\Mor(a,n\tot m)$. Then,
\begin{align*}
(\io\ot\io\ot\gamma)(X^n_{13}X^m_{23})(S\ot 1)&=(\vphi(S)\ot 1)(\io\ot\gamma)(X^a)\\&=(\vphi(S)\ot 1)(P_a \ot 1)Y^a \\ &= (\pi(S)\ot 1)Y^a = Y^{n+m}
(S \ot 1) \\ &= (Y^n(P_n \ot 1))_{13} (Y^m(P_m \ot 1))_{23} (S \ot 1) \\ &= (\io \ot \gamma)(X^n)_{13} (\io \ot \gamma)(X^m)_{23} (S \ot 1) \; .
\end{align*}
Because $(\io\ot\gamma\rho)(Y)=Y$, because $\gamma$ is multiplicative and because the coefficients of $Y$ generate $\cC$ as an algebra, it follows
that $\gamma$ is the inverse of $\rho$. So, $\rho$ is indeed a $^*$-isomorphism.
\end{proof}

\begin{remark}
A combination of Proposition 6.2.6 in \cite{bichon1} and the results in \cite{bichon2} yields an alternative proof for the fact that $A_o(F_1,F_2)
\neq 0$.
\end{remark}

\begin{remark}
Combining the co-amenability of $\SU_q(2)$ with Remark \ref{rem.coamenable}, we obtain the following. If either $F_1$ or $F_2$ is a two by two matrix
(but not necessarily both), the (unique) invariant state on the universal \cst-algebra $A_o(F_1,F_2)$ is faithful. So, the ergodic coactions of
$A_o(F_1)$ and $A_o(F_2)$ are both universal and reduced. This is somehow remarkable, because $A_o(F)$ is not co-amenable when $F$ has dimension
strictly bigger than $2$.

Moreover, still supposing that either $F_1$ or $F_2$ is a two by two matrix, we also get that $A_o(F_1,F_2)$ is a nuclear \cst-algebra. Indeed,
supposing that $F_2 \in \GL(2,\C)$, $A_o(F_1,F_2)$ is Morita equivalent with its double crossed product, i.e.\ the crossed product of the compact
operators with the co-amenable quantum group $A_o(F_2)$, yielding a nuclear \cst-algebra. See \cite{DLRZ} for details.
\end{remark}

A precise parameterisation of the unitary fiber functors (and hence, the ergodic coactions of full quantum multiplicity) on the quantum groups
$A_o(F)$, amounts to the study of matrices $F \in \GL(n,\C)$ satisfying $F \overline{F} = \pm 1$, up to the equivalence relation
\begin{equation} \label{eq.equivalence}
F_1 \sim F \quad\text{if and only if there exists a unitary $v \in
  \U(n)$ such that $F_1 = v F v^t$.}
\end{equation}
Let $F \in \GL(n,\C)$ with $F\overline{F} = \pm 1$. Denote $H=\C^n$
  and $J : H
  \recht H$ the complex conjugation. We rather look at the
  anti-linear operator $\cF = JF$, satisfying $\cF^2 = \pm 1$. In the
  case where $\cF^2 = 1$, our data come down to giving a real vector
  space together with a Hilbert space structure on the
  complexification. In the case where $\cF^2 = -1$, $H$ becomes a
  right module on the quaternions such that, restricting the
  quaternions to $\C$, we get a Hilbert space. In particular $H$ is
  even-dimensional.

It is then straightforward to provide a fundamental domain for the
equivalence relation \eqref{eq.equivalence} (see e.g.\ \cite{yamagami}). Take $F \in \GL(n,\C)$
with $F \overline{F} = \pm1$.
Let $\cF = \cJ |F|$ be
the polar decomposition of $\cF$. Then, $\cJ$ is an anti-unitary,
$\cJ^2 = \pm 1$ and $\cJ |F| \cJ^* = |F|^{-1}$. Define $H_<$ as the
subspace of $H$ spanned by the eigenvectors of $|F|$ with eigenvalue
$\lambda < 1$. Define $H_> = \cJ H_<$, which is as well the
subspace of $H$ spanned by the eigenvectors of $|F|$ with eigenvalue
$\lambda > 1$. Finally, let $H_1$ be the subset of eigenvectors of
$|F|$ with eigenvalue $1$. Take an orthonormal basis
$\xi_1,\ldots,\xi_k$ of $H_<$ of eigenvectors of $|F|$ with
eigenvalues $0 < \lambda_1 \leq \cdots \leq \lambda_k < 1$.

If $F
\overline{F} = 1$, we have $\cJ^2=1$ and we take an orthonormal basis
$\mu_1,\ldots,\mu_{n-2k}$ for $H_1$ of real vectors: $\cJ \mu_i =
\mu_i$. If $(e_i)$ denotes the standard basis for $\C^n$ and
$w : (e_i) \recht (\xi_i,\cJ \xi_i,\mu_i)$ denotes the transition
unitary, we find that
\begin{equation}\label{eq.signplus}
w^t F w = \begin{pmatrix} 0 & D(\lambda_1,\ldots,\lambda_k) & 0 \\ D(\lambda_1,\ldots,\lambda_k)^{-1} & 0 & 0 \\ 0 & 0 & 1_{n-2k} \end{pmatrix} \; .
\end{equation}
Here, $D(\lambda_1,\ldots,\lambda_k)$ denotes the diagonal matrix with
the $\lambda_i$ along the diagonal.

If $F
\overline{F} = -1$, we have $\cJ^2=-1$, $H_1$ has even dimension and
we take an orthonormal basis $\mu_1,\ldots,\mu_r,$ $\cJ \mu_1,\ldots,\cJ
\mu_r$ for $H_1$. If $w : (e_i) \recht (\xi_i,\mu_i,\cJ \xi_i,\cJ
\mu_i)$ denotes the transition unitary, we find that
\begin{equation}\label{eq.signmin}
w^t F w = \begin{pmatrix} 0 & D(\lambda_1,\ldots,\lambda_{n/2}) \\
  -D(\lambda_1,\ldots,\lambda_{n/2})^{-1} & 0 \end{pmatrix} \; ,
\end{equation}
where $0 < \lambda_1 \leq \cdots \leq \lambda_{n/2} \leq 1$.

Since the spectrum of $F^*F$ is invariant under the equivalence
relation \eqref{eq.equivalence}, a fundamental domain is given by the matrices in \eqref{eq.signplus}
with $2k \leq n$ and $0 < \lambda_1 \leq \cdots \leq \lambda_k < 1$
and the matrices in \eqref{eq.signmin} with $0 < \lambda_1 \leq \cdots \leq \lambda_{n/2} \leq 1$.

\begin{corollary} \label{cor.largemult}
Let $0 < q \leq 1$. For every even natural number $n$ with $2 \leq n \leq q + \frac{1}{q}$, the quantum group $\SU_q(2)$ admits an ergodic coaction
of full quantum multiplicity such that the multiplicity of the fundamental corepresentation is $n$. If $-1 \leq q < 0$, the same statement holds for
every natural number $n$ with $2 \leq n \leq \bigl|q + \frac{1}{q}\bigr|$.
\end{corollary}

\begin{corollary} \label{cor.classao}
Let $F \in \GL(n,\C)$ with $F \overline{F} = c 1$ and $c \in
\R$. Denote $\cG = A_o(F)$. For all $F_1 \in \GL(n_1,\C)$ satisfying
$F_1 \overline{F}_1 = c 1$ and $\Tr(F_1^* F_1) = \Tr(F^* F)$, we
denote by $\sde_{F_1}$ the coaction of $\cG$ on $A_o(F,F_1)$ defined
in Theorem \ref{thm.aoff}.
\begin{itemize}
\item Up to isomorphism, the $\sde_{F_1}$ yield all universal ergodic
  coactions of full quantum multiplicity of $\cG$. Moreover,
  $\sde_{F_1}$ is isomorphic with $\sde_{F_2}$ if and only if
  $n_1=n_2$ and there
  exists a unitary $v \in \U(n_1)$ such that $F_2 = v F_1 v^t$.
\item For all $F_1$ as above with $n_1=n$, we denote by $\Om(F_1)$
  the unitary $2$-cocycle on the dual of $\cG$ associated with the
  unitary fiber functor $\vphi_{F_1}$. The $\Om(F_1)$ describe, up to
  coboundary, all unitary $2$-cocycles on the dual of $\cG$. Moreover
  $\Om(F_1)$ and $\Om(F_2)$ differ by a coboundary if and only if
there
  exists a unitary $v \in \U(n)$ such that $F_2 = v F_1 v^t$.
\end{itemize}
\end{corollary}

\begin{corollary}
Every unitary $2$-cocycle on the dual of $\SU_q(2)$ is a coboundary.
\end{corollary}

\section{Monoidal equivalence for $A_u(F)$} \label{sec.au}

In this section, we prove, for the unitary quantum groups $A_u(F)$
studied by Banica \cite{banica1}, analogous results as for $A_o(F)$.

Again, we give a complete classification of unitary fiber functors,
monoidally equivalent quantum groups, ergodic coactions of full
quantum multiplicity and $2$-cohomology.

Recall the following definition.

\begin{definition}
For all $n\in\N$ and $F\in \GL(n,\C)$, we define $A_u(F)$ as the universal compact quantum group generated by the coefficients of the corepresentation
$U \in M_n(\C) \ot A_u(F)$ with relations
$$U \;\;\text{and}\;\; (F \ot 1)\overline{U}(F^{-1} \ot 1)
\;\;\text{are unitary}\; .$$
Then, $(A_u(F),U)$ is a compact (matrix) quantum group.
\end{definition}

Banica \cite{banica1} has shown that the irreducible corepresentations
of $A_u(F)$ can be labelled by the free monoid $\NstarN$ generated by
$\al$ and $\be$. He also computed the corresponding fusion
rules.

Defining $U^\al := U$ and $U^\be := (F \ot 1)\overline{U}(F^{-1} \ot 1)$, the spaces $\Mor(\eps,U^\al \tot U^\be)$ and $\Mor(\eps,U^\be \tot U^\al)$
are one-dimensional and generated by
$$t_F := \sum_i e_i \ot F e_i \quad\text{resp.}\quad s_F := \sum_i e_i
\ot \oF^{-1} e_i \; .$$

\begin{theorem} \label{thm.fiberau}
Let $F \in \GL(n,\C)$ be normalized such that $\Tr(F^*F) =
\Tr((F^*F)^{-1})$. Let $\cG = A_u(F)$.
\begin{itemize}
\item If $F_1 \in \GL(n_1,\C)$ satisfies $\Tr(F_1^*F_1) =
\Tr((F_1^*F_1)^{-1}) = \Tr(F^*F)$, there exists a unitary fiber
functor $\vphi_{F_1}$ on $\cG$, uniquely determined up to isomorphism,
such that
$$\vphi(t_F) = t_{F_1} \quad\text{and}\quad \vphi(s_F) = s_{F_1} \; .$$
\item Every unitary fiber functor $\vphi$ on $\cG$ is isomorphic with
  one of the form $\vphi_{F_1}$. Moreover, $\vphi_{F_1}$ is isomorphic
  with $\vphi_{F_2}$ if and only if $n_1=n_2$ and there exist unitary
  elements $v,w \in \U(n_1)$ such that $F_2 = v F_1 w$.
\end{itemize}
\end{theorem}

\begin{proof}
Let $\NstarN$ be the free monoid
generated by $\al$ and $\be$. Denote by $e$ the empty word. Elements
of $\NstarN$ are words in $\al$ and $\be$.

Take a parameter $c > 0$. Let $(A(p,q))_{p,q \in \NstarN}$ be the
universal graded \cst-algebra generated by elements
$$V_x(p,q) \in A(pq,pxq) \quad\text{for}\quad p,q \in \NstarN, \;\; x
\in  \{ \al\be,\be\al \} $$ with relations (denoting by $1_p$ the unit of the \cst-algebra $A(p) := A(p,p)$)
\begin{align*}
V_x(p,q)^* V_x(p,q) &= 1_{pq} \\
V_{\al\be}(p,\al q)^* V_{\be\al}(p\al,q) &= c \; 1_{p\al q} \\
V_{\be\al}(p,\be q)^* V_{\al\be}(p\be,q) &= c \; 1_{p\be q} \\
V_y(p,qxr) V_x(pq,r) &= V_x(pyq,r) V_y(p,qr) \\
V_x(p,qyr)^* V_y(pxq,r) &= V_y(pq,r) V_x(p,qr)^*
\end{align*}

Take $F \in \GL(n,\C)$ normalized in such a way that $\Tr(F^*F) =
\Tr((F^*F)^{-1})$. Put $c = \Tr(F^* F)$. Consider the
  quantum group $\cG=(A_u(F),\de)$. Define for every $p \in \NstarN$
  the unitary corepresentation $U^p$ of $\cG$ inductively by $U^{pq} :=
  U^p \tot U^q$.

Defining $t \in \Mor(\eps,U^{\al\be})$ and $s \in
\Mor(\eps,U^{\be\al})$ by the formulas $t= \frac{1}{\sqrt{\Tr(F^*F)}}
t_F$ and $s= \frac{1}{\sqrt{\Tr(F^*F)}} s_F$, we get a natural $^*$-homomorphism
$$\pi : (A(p,q))_{p,q \in \NstarN} \recht (\Mor(U^p,U^q))_{p,q \in \NstarN}$$
given by $\pi(V_{\al\be}(p,q)) = 1_p \ot t \ot 1_q$ and
$\pi(V_{\be\al}(p,q))=1_p \ot s \ot 1_q$. It follows from Proposition
4, the proof of Th{\'e}or{\`e}me 1 and Proposition 3 in \cite{banica1} that $\pi$ is an isomorphism of
\cst-algebras.

Take $F_1 \in \GL(n_1,\C)$ satisfying $\Tr(F_1^*F_1) =
\Tr((F_1^*F_1)^{-1}) = \Tr(F^*F)$. Write $K^\al=K^\be=\C^{n_1}$ and
define inductively $K^p$, for all $p \in \NstarN$ such that $K^{pq} =
K^p \ot K^q$. We take $K^e=\C$. From the preceding discussion, we
obtain a faithful $^*$-homomorphism
$$\pi : (\Mor(U^p,U^q))_{p,q \in \NstarN} \recht
(\B(K^p,K^q))_{p,q \in \NstarN}$$ satisfying $\pi(t) = \frac{1}{\sqrt{\Tr(F_1^*F_1)}} t_{F_1}$, $\pi(s) = \frac{1}{\sqrt{\Tr(F_1^*F_1)}} s_{F_1}$ and
$\pi(1 \ot S \ot 1) = 1 \ot \pi(S) \ot 1$ for all $S$.

We choose a concrete identification of $\cGh$ with $\NstarN$ as follows. We define, for $p \in \NstarN$, $P_p \in \Mor(U^p,U^p)$ as the unique
projection satisfying $P_p T = 0$ for all $r \in \NstarN$ with $\length r < \length p$ and all $T \in \Mor(U^r,U^p)$. We define $U_p$ as the
restriction of $U^p$ to the image of $P_p$. We then identify
$$\Mor(p_1 \tot \cdots \tot p_r, q_1 \tot \cdots \tot q_k) = (P_{q_1} \ot
\cdots \ot P_{q_k}) \Mor(U^{p_1\cdots p_r},U^{q_1 \cdots q_k})(P_{p_1} \ot \cdots \ot P_{p_r}) \; .$$ Define $H_{\vphi(p)} :=
\pi(P_p)K^p$ and define $\vphi(S)$ by restricting $\pi$ to $\Mor(p_1 \tot \cdots \tot p_r, q_1 \tot \cdots \tot q_k)$. It is obvious that $\vphi$ is
a unitary fiber functor.

The converse statement is proven in exactly the same way as in the
proof of Theorem \ref{thm.fiberao}.
\end{proof}

The next two results are proven in exactly the same way as the corresponding results for $A_o(F)$.

\begin{corollary} \label{cor.moneqau}
Let $F\in \GL(n,\C)$ with $\Tr(F^*F) = \Tr((F^*F)^{-1})$ and consider $\cG = A_o(F)$. A compact quantum group $\cG_1$ is monoidally equivalent with
$\cG$ if and only if there exists $F_1 \in \GL(n_1,\C)$ satisfying $\Tr(F_1^*F_1) = \Tr((F_1^*F_1)^{-1}) = \Tr(F^*F)$ such that $\cG_1 \cong
A_o(F_1)$.
\end{corollary}

So, we exactly know when the compact quantum groups $A_u(F_1)$ and $A_u(F_2)$ are monoidally equivalent. If this is the case, Proposition
\ref{prop.pair} provides us with a universal \cst-algebra $B_u$ and a pair of ergodic coactions of full quantum multiplicity. It is again possible to
give an intrinsic description of this \cst-algebra $B_u$. The proof is analogous to the proof of Theorem \ref{thm.aoff}. Again, the fact that
$A_u(F_1,F_2) \neq 0$ can be deduced from Proposition 6.2.6 in \cite{bichon1} and the results in \cite{bichon3}.

\begin{theorem} \label{thm.auff}
Let $F_i \in \GL(n_i,\C)$ by such that $\Tr(F_1^*F_1) =
\Tr((F_1^*F_1)^{-1}) = \Tr(F_2^*F_2) =
\Tr((F_2^*F_2)^{-1})$.
\begin{itemize}
\item Denote by $A_u(F_1,F_2)$ the universal unital \cst-algebra
  generated by the coefficients of
$$X \in M_{n_2,n_1}(\C) \ot A_u(F_1,F_2) \quad\text{with relations}\quad
X \;\;\text{and}\;\; (F_2 \ot 1) \overline{X} (F_1^{-1} \ot 1) \;\;\text{are unitary}\; .$$ Then, $A_u(F_1,F_2) \neq 0$ and there exists a unique
pair of commuting universal ergodic coactions of full quantum multiplicity, $\sde_1$ of $A_u(F_1)$ and $\sde_2$ of $A_u(F_2)$, such that
$$(\io \ot \sde_1)(X) = X_{12} (U_1)_{13} \quad\text{and}\quad (\io
\ot \sde_2)(X) = (U_2)_{12} X_{13} \; .$$
Here, $U_i$ denotes the fundamental corepresentation of $A_u(F_i)$.
\item $(A_u(F_1,F_2),\sde_1,\sde_2)$ is isomorphic with the \cst-algebra $B_u$ and the coactions thereon given by Proposition \ref{prop.pair} and the
monoidal equivalence $A_u(F_1) \meq A_u(F_2)$ of Corollary \ref{cor.moneqau}.
\end{itemize}
\end{theorem}

\begin{remark}
Exactly as in Corollary \ref{cor.classao}, a combination of Theorems
\ref{thm.fiberau} and \ref{thm.auff} gives a complete classification
of the ergodic coactions of full quantum multiplicity of $A_u(F)$ and
of the $2$-cohomology of the dual of $A_u(F)$.
\end{remark}

A precise parameterisation of the unitary fiber functors on the quantum groups $A_u(F)$ is easy. If $F_1,F \in \GL(n,\C)$, we write
$$F_1 \sim F \quad\text{if and only if there exist unitary $v,w \in \U(n)$ such that $F_1 = v F w$.}$$
We study matrices $F \in \GL(n,\C)$ satisfying $\Tr(F^* F) = \Tr((F^* F)^{-1})$ up to the equivalence relation $\sim$. It is obvious that for any
such $F$, there exist unique $0 < \lambda_1 \leq \cdots \leq \lambda_n$ satisfying $\sum_i \lambda_i^2 = \sum_i \lambda_i^{-2}$ such that $F \sim
D(\lambda_1,\ldots,\lambda_n)$. Here $D(\lambda_1,\ldots,\lambda_n)$ denotes again the diagonal matrix with the $\lambda_i$ along the diagonal.

\end{document}